\documentclass[12pt,a4paper]{amsart}
\usepackage{latexsym}
\usepackage[dvipdfmx]{graphicx}
\usepackage{epsfig} 
\usepackage{amsmath} 
\usepackage{amsfonts, amssymb} 
\usepackage{amsthm,enumerate}
\usepackage{color}
\usepackage{mathrsfs}
\usepackage{amscd}      
\usepackage{ulem}
\usepackage{euscript}
\usepackage{arydshln}
\usepackage{shuffle}
\usepackage[all]{xy}
\usepackage{here}

\usepackage[
	a4paper,
	top=25mm,
	bottom=20mm,
	left=25mm,
	right=25mm,
]{geometry}

\theoremstyle{definition}
\newtheorem{Def}{Definition}[section]

\newtheorem{Eg}{Example}[section]
\newtheorem{Rm}{Remark}[section]

\theoremstyle{plain}
\newtheorem{Prop}[Def]{Proposition}

\newtheorem{Thm}[Def]{Theorem}

\numberwithin{equation}{section}

\newcommand{\authorfootnotesA}{\renewcommand\thefootnote{$\flat$}}%
\newcommand{\authorfootnotesB}{\renewcommand\thefootnote{$\sharp$}}%
\newcommand{\authorfootnotesC}{\renewcommand\thefootnote{$\diamond$}}%
\newcommand{\authorfootnotesD}{\renewcommand\thefootnote{$\S$}}%

\begin{document}

\begin{center}
	\LARGE 
	Simulation Example of a Black Noise \par \bigskip

	\normalsize
	\authorfootnotesA
	T.~Amaba\footnote{%
	This work was partially supported by funding from Fukuoka University (Grant No. 197102) and by JSPS KAKENHI Grant Number 22K03345.%
	}\textsuperscript{1},
	\authorfootnotesB
	T.~Aoyama\textsuperscript{2},
	\authorfootnotesC
	S.~Araki\textsuperscript{3}
	and
	\authorfootnotesD
	S.~Eguchi\textsuperscript{4}

	\textsuperscript{1}Fukuoka University,
	8-19-1 Nanakuma, J\^onan-ku, Fukuoka, 814-0180, Japan.\par \bigskip

	\textsuperscript{2}Okayama University of Science, 1-1 Ridai-cho, Kita-ku, Okayama, 700-0005, Japan.\par \bigskip
	
	\textsuperscript{3}Fukuoka University Ohori Senior High School, 1-12-1 Ropponmatsu, Chu\^o-ku, Fukuoka, 810-0044, Japan.\par \bigskip
	
	\textsuperscript{4}Aomori University, 2-3-1 K\^obata, Aomori, Aomori, 030-0943, Japan.\par \bigskip

	
	\email{(T.~Amaba) fmamaba@fukuoka-u.ac.jp}
	\email{(T.~Aoyama) aoyama@ous.ac.jp}
	\email{(S.~Araki) s.araki.ho@adm.fukuoka-u.ac.jp}
	\email{(S.~Eguchi) eguchi.math@gmail.com}

	\today
\end{center}

\begin{quote}{\small {\bf Abstract.}
B. Tsirelson and A. M. Vershik (1998) introduced the notion of a mathematical noise, which possesses completely opposite properties to those of a white noise. Afterward, B. Tsirelson (2004) called this noise: `black noise.' In this paper, we provide a method to simulate black noise using a modified Bayesian convolutional neural network. Then we study the behavior of black noise both numerically and visually.
}
\end{quote}

\renewcommand{\thefootnote}{\fnsymbol{footnote}}



\section{Introduction}
\label{Intro}
In Error Theory, (Gaussian) white noise has been used for a long time as a mathematical model of fluctuations caused by a number of random errors in various applied fields such as communication theory, informatics, mathematical finance, statistical physics, and theoretical physics among others.

In particular, in audio engineering, electronics, and physics, white noise is represented through a flat spectrum in the frequency domain. Colored noises are defined by non-flat spectrums. More specifically,

{\sc White Noise:}
If random amplitudes of frequencies appear the spectrum which can be considered to be independent and identically distributed according to a mean-zero Gaussian distribution, the corresponding noise is called a Gaussian {\it white noise}. Since this noise has frequency spectrums in the range that human beings may perceive, we can recognize that it oscillates randomly as a wave after taking Fourier inversion of the spectral series.
Gaussian white noise is particularly compatible with $L_{2}$ theory, and the work of
It\^{o}
(\cite{Ito1,Ito2}),
Kunita–Watanabe
(\cite{KuWa}),
Ogawa
(\cite{Oga1,Oga2}),
and
Watanabe
(\cite{Wa})
is widely known in the context of $L_{2}$ theory of stochastic analysis.

{\sc Black Noise:}
In the case of sound, one says that it is `{\it black noise}' if the amplitudes of frequencies are non-zero only outside of the range that human ears may perceive. That is, we cannot hear this noise although they exist physically. In this article, we modify this definition as an object whose finite frequency spectrum values are zero even though it does have nonzero energy! As a result, the black noise oscillates too fast to be recognized by standard measurements such as a spectrum oscillator or, equivalently to say, any linear sensors.

A linear sensor is the aggregation of a signal which mathematically can be characterized as a random field $\omega (x)$ taking values in space or space-time.
A sensor associates to a field $\omega (x)$ the value $\int f(x) \omega (x) \mathrm{d}x$ where $f$ belongs to a class of function which depends on the context.
For example, we may take
\begin{itemize}
\item[$\circ$]
$f(x) = \exp (-\sqrt{-1}\, \xi x )$
when we want to know the value of the spectrum of $\omega$ at the frequency $\xi$,

\item[$\circ$]
$f(x) = \delta_{0} (x)$
when we want to detect the value $\omega (0)$.

\end{itemize}

Any electromagnetic field is sensed using the above transformation at a finite frequency spectrum and then the total amount of the sensed frequency spectrum gives the energy, and eventually, the field itself is identified by data obtained this way.

It is said that hypothetical dark matter does not interact with electromagnetic fields. This should imply that it cannot be detected by any linear sensors despite their nonzero energy. Black noise have features similar to this requirement. Actually, in \cite{Ka2}, G.~Kalai (2007) speculated that `dark energy is a black noise.' (See also \cite{Ka1}.)
In \cite[p.~242]{Tsi3}, Tsirelson quotes a comment of Shnirelman (\cite{Shn}) on the paradoxical motion of an ideal incompressible fluid as one of the motivations for his discussion of black noise: `... very strong external forces are present, but they are infinitely fast oscillating in space and therefore are indistinguishable from zero in the sense of distributions. The smooth test functions are not ``sensitive" enough to ``feel" these forces.' Based on this, he states that a fluid could be a nonlinear sensor that senses black noises.

Black noise was first constructed by Tsirelson–Vershik in \cite{TsiVer}. Afterwards, some examples of similar kinds of such noise, they also call them black, were given by Le~Jan–Raimond (\cite{LeRa1,LeRa2}), Warren–Watanabe (\cite{WaWa}) and Watanabe (\cite{Wa2}).
They all followed the techniques of stochastic flows to give them, which are quite different from the one given in \cite{TsiVer}.
The term `black noise' first appeared in \cite{Tsi1}.
For any Markov-diffusion-type stochastic differential equation, in which a unique strong solution exists for any initial point, the associated stochastic flow determines a noise through a certain standard procedure and this noise is known to be white (in the sense of Tsirelson, see Example~\ref{Eg:white_noise} below). 
Therefore black noises do not appear in the nondegenerate setting that is usually common when dealing with stochastic differential equations. 
The theory of the part of the noise described by stochastic flows is very nicely organized, but it seems to be difficult to give a natural interpretation as some existing fields such as scalar fields to noises constructed by following the techniques of stochastic flows.

Among the noises associated with stochastic flows, the most famous black noise is the system of coalescing Brownian motions, also known as {\it Arratia flow} (also called {\it Brownian web}). 
These are one-dimensional Brownian motions starting from each point in the real line, moving independently of each other until they collide, and then they both begin to move together after the collision. 
The dynamics can be described by stochastic differential equations (due to the paper Warren–Watanabe~\cite{WaWa}), but the diffusion coefficients degenerate violently.
We refer the reader to \cite{AIW} for a review of other noises associated with stochastic flows.

It is also known that any scaling limit of the critical planar percolation defines a (two-dimensional version of) black noise. This was conjectured by Tsirelson and affirmatively resolved by Schramm–Smirnov–Garban (\cite{SSG}).
We may find very few afterward examples of black noises such as the result of rethinking Arratia flow as a noise in a two-dimensional time-space (see Ellis–Feldheim \cite{ElFe}). 
Still, examples of black noises (in particular, in two-dimensional cases) are not abundant.

Though it is easy to recognize the shape of the white noise, it is difficult to describe how a black noise behaves.
As are shown in some examples above, black noise is a field that appears in very singular systems.
In the case of two-dimensional (generally multidimensional) systems, as suggested by the percolation example, it is related not only to a singular system but also to some symmetry of the system such as conformal covariance or renormalization group of the system.

The above background motivated us to develop a method for numerical experiments using black noise which is known to be ultimately harsh in the sense that, as we mentioned above, they are unobservable with all linear detectors such as the Fourier transformation. Although the black noise we treat in this article may be regarded as a random field, their realizations are no longer functions, tensor fields, or even distributions. These realizations will be described as some elements of a net constructed by function spaces following Tsirelson and Vershik~\cite{TsiVer}.

In this article, we implement the black noise constructed by Tsirelson and Vershik (\cite{TsiVer}), give a displayable approximating sample of them, and claim the blackness by checking a part of the equivalent condition of blackness due to Tsirelson (2004).
For this purpose, we will use Proposition~\ref{Prop:the_contr} to be shown later. 
This proposition contributes technically in the current paper as a clue to construct a Bayesian statistical model by following the idea of Bayesian neural network.
To our best knowledge, this is the first study that tries to simulate black noise and describe its behavior.

\section{White Noise and Black Noise}

\subsection{Definition of noise, white noise and black noise}

We first introduce the notion of noise by following Tsirelson~\cite{Tsi2} and Akahori et al.~\cite{AIW}.

\begin{Def} 
For a complete probability space
$( \Omega , \mathcal{F}, \mathbf{P} )$,
a family of complete sub-$\sigma$-fields
$\{ \mathcal{F}_{s,t} \}_{-\infty < s \leqslant t < \infty}$
of $\mathcal{F}$
with
$
\mathcal{F}
= \mathcal{F}_{-\infty , \infty}
:= \vee_{s \leqslant t} \mathcal{F}_{s,t}
$
is called a {\it noise} if
\begin{itemize}
\item[(1)]
for all $s \leqslant t \leqslant u$, it holds that
$
\mathcal{F}_{s,t} \otimes \mathcal{F}_{t,u} = \mathcal{F}_{s,u}
$
by which we mean
$
\mathcal{F}_{s,t} \vee \mathcal{F}_{t,u} = \mathcal{F}_{s,u}
$
and
$\mathcal{F}_{s,t}$ and $\mathcal{F}_{t,u}$ are independent,

\item[(2)]
there exists a one-parameter measurable group
$\{ T_{h} \}_{h \in \mathbb{R}}$,
which consists of
$
\mathcal{F}_{-\infty , \infty}/\mathcal{F}_{-\infty , \infty}
$-measurable maps
$T_{h}: \Omega \to \Omega$,
such that
$T_{h}$ is bijective, preserves $\mathbf{P}$,
and satisfies
$
T_{h} ( \mathcal{F}_{s,t} ) = \mathcal{F}_{s+h,t+h}
$
for all $h \in \mathbb{R}$.
\end{itemize}
We denote this noise by
$
\mathbf{N}
=
(
	\Omega ,
	\{ \mathcal{F}_{s,t} \}_{-\infty <s \leqslant t < \infty},
	\mathbf{P},
	\{ T_{h} \}_{h \in \mathbb{R}}
)
$
or
$
\mathbf{N}
=
(
	\{ \mathcal{F}_{s,t} \}_{s \leqslant t},
	\{ T_{h} \}_{h \in \mathbb{R}}
)
$
as a shorthand notation.
\end{Def} 

In other words, a noise is a family of $\sigma$-fields parametrized by time intervals with Markov property in the sense of (1) and covariance regarding parallel shifts along the time axis.

\begin{Eg}[White Noise] 
\label{Eg:white_noise} 

Let $W = (w(t))_{t \in \mathbb{R}}$ be a $d$-dimensional
Wiener process on
$\Omega = C( \mathbb{R} \to \mathbb{R}^{d} )$
(Recall that $W_{t} (w) = w(t)$).
Set
$
\mathcal{F}_{s,t}^{W}
:=
\sigma [ W_{v} - W_{u}: s \leqslant u \leqslant v \leqslant t ]
$
and
$T_{h}: \Omega \ni w \mapsto w ( \bullet - h ) \in \Omega$.
Then
$
\mathbf{N}^{W}
=
(
	\{ \mathcal{F}_{s,t}^{W} \}_{s \leqslant t},
	\{ T_{h} \}_{h \in \mathbb{R}}
)
$
forms a noise called a
{\it $d$-dimensional white noise}
(in the sense of Tsirelson).
Note that $\omega (t) = w^{\prime} (t)$, $t \in \mathbb{R}$
(here, $w^{\prime} (t)$ is the distributional derivative of $w(t)$ with respect to $t$),
is the Gaussian white noise,
which is widely used in various scientific fields.

\end{Eg} 

We have another well-known process which describes certain randomly moving or appearing particles and is regarded as quite different from white noise above.

\begin{Eg}[Poisson Noise] 

Let $p$ be a stationary Poisson point process\footnote{
A pair $(p, \mathbf{D}_{p})$ is called a
{\it point function}
(which we denote by $p$ again)
on a measurable space $\mathbf{X}$
if
$\mathbf{D}_{p} \subset \mathbb{R}$
is countable and
$p: \mathbf{D}_{p} \to \mathbf{X}$.
We denote by
$\Pi_{\mathbf{X}}$
the set of all point functions on $\mathbf{X}$.
Then the mapping
$
\Pi_{\mathbf{X}} \ni p \mapsto N_{p} ( \mathrm{d}t \mathrm{d}x )
:=
\# \{ \mathrm{d}t \cap \mathbf{D}_{p}: p(x) \in \mathrm{d}x \}
$
induces a measurable structure on $\Pi_{\mathbf{X}}$.
A $\Pi_{\mathbf{X}}$-valued random variable $p$
is called a
{\it stationary Poisson point process}
on a $\sigma$-finite measure space $(\mathbf{X}, n)$
with the
{\it characteristic measure}
$n$ if
$
N_{p} ( \mathrm{d}t \mathrm{d}x )
$
is a Poisson random measure with intensity
$
\mathrm{d}t \otimes n ( \mathrm{d}x )
$.
}
on a Polish space $\mathbf{X}$,
and
set up on the canonical space
\begin{equation*}
\Omega := \Pi_{\mathbf{X}}
=
\{ \text{point function $(p, \mathbf{D}_{p})$ on $\mathbf{X}$} \} .
\end{equation*}
Then by setting
$
\mathcal{F}_{s,t}^{p}
:=
\sigma [ p(u): s \leqslant u \leqslant t ]
$
and
$
T_{h}: \Pi_{\mathbf{X}} \ni (p, \mathbf{D}_{p})
\mapsto
( p(\bullet -h), \mathbf{D}_{(T_{h})^{*}p} ) \in \Pi_{\mathbf{X}}
$,
where
$
\mathbf{D}_{(T_{h})^{*}p}
:=
\{ t : t-h \in \mathbf{D}_{p} \}
$,
we see that
$
(T_{h})_{*}p
:= (T_{-h})^{*}p
= p ( \bullet + h )
$
almost surely
and
$
\mathbf{N}^{p}
=
(
	\{ \mathcal{F}_{s,t}^{p} \}_{s \leqslant t},
	\{ T_{h} \}_{h \in \mathbb{R}}
)
$
is a noise.
This is called a {\it Poisson noise}
(or a {\it shot noise}) which occurs in photon counting in optical devices and so on.
\end{Eg} 

One may construct a noise using a combination of white and Poisson noises as below.

\begin{Eg}[Classical Noise] 
\label{Eg:Levy} 
Any arbitrary L\'{e}vy process
$X = (X_{t})_{t \in \mathbb{R}}$
defines a noise
$
\mathbf{N}
=
(
	\{ \mathcal{F}_{s,t} \}_{s \leqslant t},
	\{ T_{h} \}_{h \in \mathbb{R}}
)
$
by setting
$
\mathcal{F}_{s,t}
=
\sigma ( X_{v}-X_{u} : s \leqslant u \leqslant v \leqslant t )
$.
The transformation $T_{h}$ is defined similarly to that in Example~\ref{Eg:white_noise}.
The L\'{e}vy-It\^{o} decomposition theorem implies that
there exist a white noise $\mathbf{N}^{W}$
and a Poisson noise $\mathbf{N}^{p}$
such that
$
\mathcal{F}_{s,t}
=
\mathcal{F}_{s,t}^{W} \otimes \mathcal{F}_{s,t}^{p}
$.
Noises obtained in this way are called
{\it classical noises}.
\end{Eg} 

\begin{Def} 

Let
$
\mathbf{N}
=
( \{ \mathcal{F}_{s,t} \}_{s\leqslant t}, \mathbf{P}, \{ T_{h} \}_{h \in \mathbb{R}} )
$
be a noise.
\begin{itemize}
\item[(1)]
A function
$f \in L_{2}(\mathcal{F}_{-\infty , \infty})$
is called a
{\it first order chaos}
(or {\it $\mathbf{P}$-integral})
if
for any
$s \leqslant t \leqslant u$,
it holds that
$
\mathbf{E} [ f \vert \mathcal{F}_{s,u} ]
=
\mathbf{E} [ f \vert \mathcal{F}_{s,t} ]
+
\mathbf{E} [ f \vert \mathcal{F}_{t,u} ]
$.

\item[(2)]
The
noise $\mathbf{N}$ is called
{\it black}
if it does not have any first order chaos except $0 \in L_{2}( \mathcal{F}_{-\infty , \infty} )$.
\end{itemize}
\end{Def} 

Roughly speaking, first order chaos is an axiomatization of the Wiener integrals as shown in the following example.

\begin{Eg}[First order chaos as linear sensors]
\label{Eg:1stChaos} 

In the case of white noise,
$\mathbf{N}^{W}$,
the space of first chaos random variables consists of
all Wiener integrals
$
\int_{-\infty}^{\infty} h(x) \mathrm{d} w(x)
$,
$h \in L_{2} ( \mathbb{R} )$.
The formal derivative
$
\omega (x) = w^{\prime} (x)
$
is what originally is called a white noise in many fields.
By using this notation, the previous integral can be written as
$
\int_{-\infty}^{\infty} h(x) \omega (x) \mathrm{d} x
$
which is one type of linear sensor as described in the introduction
(Section~\ref{Intro}).
\end{Eg} 

In the above example, the first chaos can be regarded as a generalization of linear sensors. In this sense, we may say that our goal is to find a setting where all first chaos Wiener integrals will be zero. This will imply that black noise cannot be detected by any linear sensors.

A subset $A \subset \mathbb{R}$ is called an
{\it elementary set}
if $A = \cup_{i=1}^{m} [ s_{i}, t_{i} ]$
for some $m \in \mathbb{N}$ and
$
-\infty < s_{1} < t_{1} < \cdots < s_{m} < t_{m} < \infty
$.
We denote by
$\EuScript{C}$
the set of all compact sets in $\mathbb{R}$.
(By convention, we
allow
$\EuScript{C}$ to include the empty set $\varnothing$.)
By equipping the space $\EuScript{C}$ with the Hausdorff distance, it becomes a Polish space, and we denote by $\mathcal{B} ( \EuScript{C} )$ the Borel $ \sigma $-field on $ \EuScript{C} $ equipped with the Hausdorff-distance topology.
We set
$
\EuScript{C}_{A}
:=
\{ M \in \EuScript{C}: M \subset A \}
$.

The following result is due to Tsirelson (2004).

\begin{Prop}[{\cite[Theorem~3d12]{Tsi2}}] 
\label{Prop:subsp} 

For any noise
$
\mathbf{N}
=
( \{ \mathcal{F} \}_{s\leqslant t}, \{ T_{h} \}_{h \in \mathbb{R}} )
$,
there is a unique mapping
\begin{equation*}
\mathcal{B} ( \EuScript{C} )
\ni \EuScript{M}
\mapsto
(\text{a closed subspace $H_{\EuScript{M}}$ of
$
H := L_{2}( \Omega , \mathcal{F}_{-\infty , \infty}, \mathbf{P} )
$})
\end{equation*}
such that
\begin{itemize}
\item[{\rm (i)}]
for any set $A=\cup_{i=1}^{m} [ s_{i}, t_{i} ]$ with $
-\infty < s_{1} < t_{1} < \cdots < s_{m} < t_{m} < \infty
$, we have 
$
H_{\EuScript{C}_{A}}
=
L_{2} ( \mathbb{F} (A) )
$
where
$
\mathbb{F} (A)
:=
\mathcal{F}_{s_{1},t_{1}}
\otimes \cdots
\otimes \mathcal{F}_{s_{m},t_{m}}.
$

\item[{\rm (ii)}]
if
$
\EuScript{M}_{1}, \EuScript{M}_{2},
\ldots
\in \EuScript{C}
$
are pairwise disjoint,
then
$
H_{\EuScript{M}_{1} \cup \EuScript{M}_{2} \cup \cdots}
=
H_{\EuScript{M}_{1}} \oplus H_{\EuScript{M}_{2}} \oplus \cdots
$.
\end{itemize}
\end{Prop} 

For each $\EuScript{M} \in \mathcal{B}(\EuScript{C})$,
denote by
$
P_{\EuScript{M}}
$
the orthogonal projection
from $H$ onto $H_{\EuScript{M}}$.
Then the following holds:
(1)
for each elementary set $A$ and $f \in H$,
it holds that
$
P_{\EuScript{C}_{A}} f
=
\mathbf{E} [ f \vert \mathbb{F}(A) ]
$.
(2)
For each $f \in L_{2} ( \mathcal{F}_{-\infty , \infty} )$,
the mapping
$
\mu_{f} : \mathcal{B}(\EuScript{C}) \ni \EuScript{M}
\to
\Vert P_{\EuScript{M}} f \Vert_{H}^{2}
$
defines a measure with total
variation
$
\Vert f \Vert_{H}^{2}
$.
With these notations, we have the following definition.

\begin{Def} 

The measure
$
\mu_{f}
$
is called the
{\it spectral measure}
of the noise
$
\mathbf{N}
=
( \{ \mathcal{F} \}_{s\leqslant t}, \{ T_{h} \}_{h \in \mathbb{R}} )
$
with respect to
$
f \in H = L_{2}( \Omega , \mathcal{F}_{-\infty , \infty}, \mathbf{P} )
$.

\end{Def} 

It is clear that
$
H_{\{ \varnothing \}}
=
(\text{constant random variables})
\cong
\mathbb{R}
$.
Using Proposition~\ref{Prop:subsp}, we have that for each $n \in \mathbb{N}$,
$
\EuScript{C}_{n} = \{ M \in \EuScript{C}: \# M = n \}
\in \mathcal{B} ( \EuScript{C} )$
is associated with the closed subspace
$H_{\EuScript{C}_{n}}$
of
$
H = L_{2}( \mathcal{F}_{-\infty , \infty} )
$.
Now by setting
$
\EuScript{C}_{\mathrm{finite}}
=
\{ M \in \EuScript{C}: \# M <  \infty \}
$,
Proposition~\ref{Prop:subsp}--(ii)
implies that
\begin{equation*}
H_{\EuScript{C}_{\mathrm{finite}}}
=
\mathbb{R}
\oplus H_{\EuScript{C}_{1}}
\oplus H_{\EuScript{C}_{2}}
\oplus \cdots
\end{equation*}
is an orthogonal decomposition of $H_{\EuScript{C}_{\mathrm{finite}}}$.
For a subset $K$ of $H$, we denote by $\sigma (K)$ the $\sigma$-field generated by $K$.

\begin{Prop}[{\cite[Section~6a, Theorem~6a3]{Tsi2}}] 

The subspace
$H_{\EuScript{C}_{1}}$
coincides with the space of all first order chaos,
and it holds that
$
\sigma ( H_{\EuScript{C}_{\mathrm{finite}}} )
=
\sigma ( H_{\EuScript{C}_{1}} )
$.

\end{Prop} 

\begin{Def} 

The subspace
$H_{\EuScript{C}_{n}}$
is called
the
space of
{\it $n$-th order chaos}.

\end{Def} 

With the above definition, we see that Proposition~\ref{Prop:subsp} is a manifestation of the stratification of randomness property of the noise.
In the field of Malliavin Calculus, the Malliavin derivative operator sends $H_{\EuScript{C}_{n}}$ to $H_{\EuScript{C}_{n-1}}$, namely the operator reduces the randomness at least of order one, or in the setting of quantum mechanics, the Malliavin derivative operator acts as an annihilation and its adjoint operator acts as a creation. 
They satisfy the Heisenberg commutation relation, making $H_{\EuScript{C}_{\mathrm{finite}}}$ an irreducible representation space of Heisenberg algebra.

In the case of black noise, it holds by definition that $H_{\EuScript{C}_{1}} = \{ 0 \}$, so that $H_{\EuScript{C}_{n}} = \{ 0 \}$ for all $n \in \mathbb{N}$, which implies that $H = \mathbb{R} \oplus H_{\EuScript{C} \setminus \EuScript{C}_{\mathrm{finite}}}$. 
Hence the set of $f \in H$, for which its mean is zero but the squared $L_2$-norm (also called energy) is non-zero is a subset of $H_{\EuScript{C} \setminus \EuScript{C}_{\mathrm{finite}}}$.

\subsection{Tsirelson--Vershik's construction of a black noise}
\label{Sec:Construction}

To construct a black noise,
Tsirelson and Vershik
introduced in \cite{TsiVer}
the following projective system:
\begin{equation}
\label{eq:proj_sys} 
\Omega_{0} (\mathbb{R}) \stackrel{R_{0,1}}{\leftarrow}
\Omega_{1} (\mathbb{R}) \stackrel{R_{1,2}}{\leftarrow}
\cdots \stackrel{R_{k-1,k}}{\leftarrow}
\Omega_{k} (\mathbb{R}) \stackrel{R_{k,k+1}}{\leftarrow}
\cdots ,
\end{equation}
where for every $k$,
\begin{itemize}
\item[(1)]
$
\Omega_{k} (\mathbb{R}) := C(\mathbb{R})
= (\text{the space of continuous functions on $\mathbb{R}$})
$,

\item[(2)]
$
r_{0}
>0
$, $M>1$ and
$
r_{k+1}
:=
\frac{r_k}{M}
= \frac{r_1}{M^k}
$,

\item[(3)]
$
s_{k}
:= \sum_{l=1}^{\infty} r_{k+l}
= \frac{1}{M-1} r_{k}
$,

\item[(4)]
$
R_{k-1,k}:
\Omega_{k} (\mathbb{R}) \ni \omega_k
\mapsto
\omega_{k-1} \in \Omega_{k-1}(\mathbb{R})
$
is defined by
\begin{equation}
\label{Def:R-map} 
\omega_{k-1}(x)
:=
\varphi
\left(
	L \frac{\sqrt{M-1}}{2r_k}
	\int_{x-r_k}^{x+r_k} \omega_k (y) \mathrm{d}y
\right),
\quad
x \in \mathbb{R},
\end{equation}
\end{itemize}
and where
\begin{itemize}
\item[(a)]
$L, M>0$
are sufficiently large constants
which will be discussed
in Theorem~\ref{Thm:Tsi-Ver} below,

\item[(b)]
the function
$\varphi : \mathbb{R} \to \mathbb{R}$
is fixed and we assume that it satisfies the following properties:
	\begin{itemize}
	\item[(i)]
	(Antisymmetry)
	$\varphi (u) \equiv -\varphi (-u)$,

	\item[(ii)]
	(Lipschitz condition)
	there exists a constant $C>0$ such that
	$$
	\vert \varphi (u) - \varphi (v) \vert
	\leqslant
	C \vert u-v \vert,
	\quad
	u,v \in \mathbb{R},
	$$

	\item[(iii)]
	(Nonlinearity)
	$\varphi (u) = 1$
	for $u \geqslant 1$.
	\end{itemize}
\end{itemize}

We denote the projective limit of the diagram
(\ref{eq:proj_sys})
by
\begin{equation*}
\begin{split}
&
\Omega (\mathbb{R})
:=
\varprojlim \Omega_{k} (\mathbb{R}) \\
&=
\left\{
	\omega = ( \omega_{k} )_{k=0}^{\infty}
	\in \prod_{k=0}^{\infty} \Omega_{k} (\mathbb{R})
	:
	\text{$\omega_{k-1} = R_{k-1,k} \omega_{k}$ for every $k \geqslant 1$}
\right\}
\end{split}
\end{equation*}
with
$R_{k,\infty} : \Omega (\mathbb{R}) \to \Omega_{k} (\mathbb{R})$
being the $k$-th coordinate projection from
$
\prod_{k=0}^{\infty} \Omega_{k}(\mathbb{R})
$.
We denote by $\mathbb{F}(\mathbb{R})$
the Borel-$\sigma$-field on $\Omega (\mathbb{R})$.
For every open set $U \subset \mathbb{R}$,
we endow $\Omega (\mathbb{R})$ an equivalence relation defined by
\begin{equation*}
\begin{split}
&
\omega^{\prime} = ( \omega_{k}^{\prime} )_{k=0}^{\infty}
\stackrel{U}{\sim}
\omega^{\prime\prime} = ( \omega_{k}^{\prime\prime} )_{k=0}^{\infty} \\
&\quad\quad\Longleftrightarrow
\left(\begin{array}{c}
\text{for each $k=0,1,2, \ldots$, \phantom{if $U_{-s_{k}} \neq \varnothing$ then}} \\
\text{$
U_{-s_{k}} \neq \varnothing
\Rightarrow
\omega_{k}^{\prime}\vert_{U_{-s_{k}}}
=
\omega_{k}^{\prime\prime}\vert_{U_{-s_{k}}}
$,
}
\end{array}
\right)
\end{split}
\end{equation*}
where for $r>0$,
$
U_{-r}
=
\{ x \in \mathbb{R}: \inf_{y \in \mathbb{R} \setminus U}\vert x - y \vert < r \}
$.
The quotient set
$\Omega (U) := \Omega (\mathbb{R}) / \stackrel{U}{\sim}$
inherits a measurable structure from $\Omega (\mathbb{R})$
through the natural projection
$
\rho_{U,\mathbb{R}} : \Omega (\mathbb{R}) \to \Omega (U)
$.
We denote by
$
\mathbb{F}(U) := \sigma ( \rho_{U,\mathbb{R}} )
$
and now
we have a family
$
\{ \mathbb{F} (U) \}_{\text{$U$: open in $\mathbb{R}$}}
$
of sub-$\sigma$-fields of $\mathbb{F} (\mathbb{R})$.
For each $h \in \mathbb{R}$, we define
$T_{h} : \Omega (\mathbb{R}) \to \Omega (\mathbb{R})$
by the component-wise translation by $-h$.

\begin{Rm} 
\label{Rm:CNN} 
If we discretize  (\ref{Def:R-map}) as
\begin{equation*}
\omega_{k-1}(x)
\approx
\varphi
\left(
	L \frac{\sqrt{M-1}}{2r_k}
	\sum_{\substack{
	j: \\
	x-r_{k} \leqslant y_{j} \leqslant x+r_{k}
	}}
	\omega_{k} ( y_{j} )
	\Delta y_{j}
\right) ,
\end{equation*}
the above formula looks like the basic iteration in a Convolutional Neural Network
(CNN)
by considering
\begin{itemize}
\item[$\bullet$]
$\varphi$ as the activation function,

\item[$\bullet$]
$\omega_{k}$ as the input into $k$-th layer,

\item[$\bullet$]
$L \frac{\sqrt{M-1}}{2r_k} \Delta y_{j}$
as the weights at the unit $j$.
In the context of Bayesian CNN, these weights are replaced by unit- and layer-wise independent random variables as prior distributions. 
The most commonly used distributions for the weights are Gaussian.
By taking Brownian motions
$W^{k} = ( W^{k}(y) )_{y \in \mathbb{R}}$
which are independent with respect to $ k $,
one may propose to replace  (\ref{Def:R-map}) by
\begin{equation*}
\omega_{k-1}(x)
=
\varphi
\left(
	L \frac{\sqrt{M-1}}{2r_k}
	\int_{x-r_{k}}^{x+r_{k}}
	\omega_{k} (y) \mathrm{d}W^{k}(y)
\right),\ x \in \mathbb{R}.
\end{equation*}
This proposal can be interpreted as a continuum version of a statistical machine-learning model as illustrated in Neal~\cite{Ne} and Lee et al.~\cite{LBNSPS}.
\end{itemize}
\end{Rm} 

\begin{Def}[The space $\mathcal{M}$ and $\gamma \in \mathcal{M}$] 
\label{Def:Space_M} 
\phantom{jhsdckgakyg}
\begin{itemize}
\item[(1)]
We denote by
$\mathcal{M}$
the set of probability measures $\mu$ on
$C(\mathbb{R})$
which satisfy the following properties:
	\begin{itemize}
	\item[$\circ$]
	$\mu$ is invariant under the transformations
	$
	C(\mathbb{R} )
	\ni w \mapsto w (\bullet -t) \in
	C(\mathbb{R})
	$
	for $t \in \mathbb{R}$,

	\item[$\circ$]
	$\mu$ is invariant under the transformation
	$
	C(\mathbb{R} )
	\ni w \mapsto -w \in
	C(\mathbb{R})
	$,

	\item[$\circ$]
	$\mu$ is $2$-dependent:
	for any open sets $U$ and $V$ in $\mathbb{R}$ which satisfy 
	$\inf_{(x,y) \in U \times V} \vert x-y \vert \geqslant 2$,
	$\mathbb{F}(U)$
	and
	$\mathbb{F}(V)$
	are independent under $\mu$.
	\end{itemize}

We endow on $\mathcal{M}$ the weak topology and treat it as a topological space.

\item[(2)]
We denote by
$\gamma \in \mathcal{M}$
the Gaussian measure on $C(\mathbb{R})$
of which the mean function is zero
and
the variance function is given by,
for $x,y \in \mathbb{R}$,
$$
\int_{C(\mathbb{R})}\xi (x) \xi (y) \gamma (d\xi )
=
2 \frac{\pi -2}{\pi}
\max
\big\{
	0,
	1 - \frac{1}{2}
	\vert x-y\vert
\big\} .
$$
Remark that the right-hand side of the above equality is a positive definite function of $x$ and $y$, and thus the existence of such measure is guaranteed by the Moore--Aronszajn theorem
(See e.g. \cite{Ar}).
\end{itemize}
\end{Def} 

We remark here that one may consider a $r$-dependent condition instead of a $2$-dependent condition in Definition~\ref{Def:Space_M}--(1). With some trivial changes, similar conclusions to the ones obtained here follow.

\begin{Def} 
\label{Def:SB} 
\begin{itemize}
\item[(1)]
We define
$S^{L,M} : C(\mathbb{R}) \to C(\mathbb{R})$
by
\begin{equation*}
\begin{split}
(S^{L,M} \xi)(x)
&:= 
\displaystyle \frac{1}{\sqrt{M}} \int_0^{Mx} \varphi (L\xi (y)) \mathrm{d}y,
\quad
x \in \mathbb{R}.
\end{split}
\end{equation*}

\item[(2)]
For $\mu \in \mathcal{M}$ and $L>0$,
$$
B_{L} (\mu)
:=
\int_{0}^{2} \mathbf{E}_{\mu} 
[\varphi (L\xi (0))\varphi (L\xi (x))] \mathrm{d}x.
$$
\end{itemize}

\end{Def} 

For a given measure space
$
(X , \mu)
$,
a measurable space $Y$
and a measurable mapping $f: X \to Y$,
we denote by
$f(\mu )$
the image measure of $\mu$ by $f$:
That is, the measure on $Y$ defined by
$
( f(\mu ) )( \mathrm{d}y )
=
\mu ( f \in \mathrm{d}y )
$.

The following result is due to Tsirelson and Vershik (1998).

\begin{Thm}[{\cite[Corollary~5.3, Lemma~C2, Lemma~C4]{TsiVer}}]
\label{Thm:Tsi-Ver} 
\ \par
\begin{itemize}
\item[{\rm (1)}]
For any open neighborhood in the weak topology
$\mathcal{U} \subset \mathcal{M}$
of $\gamma$,
there exist
$L,M >0$
and a probability measure $\mu$ on
$\Omega  (\mathbb{R}) $
with the following properties:
\begin{itemize}
	\item[{\rm (i)}]
	The tuple $( \Omega (\mathbb{R}) ,\{ \mathbb{F}(s,t) \}_{s \leqslant t}, \mu , \{ T_{h} \}_{h \in \mathbb{R}} )$ is a black noise.
	
	
	\item[(ii)]
	Let $\omega$ be the identity map
	(which is also called the canonical process)
	on $\Omega (\mathbb{R})$.
	Then for $k=0,1,2, \ldots $,
	the law of $\xi_{k} = ( \xi_{k} (x) )_{x \in \mathbb{R}}$
	belongs to $\mathcal{U}$,
	where $\xi_{k} (x)$ is defined by
	\begin{equation}
	\label{Def:Xi-map} 
		\begin{split}
			\xi_{k} (x)
			&:= 
			\displaystyle \frac{\sqrt{M-1}}{2r_{k+1}}
			\int_{s_k x-r_{k+1}}^{s_k x + r_{k+1}} (R_{k+1,\infty} \omega )(y)\mathrm{d}y.
		\end{split}
	\end{equation}
\end{itemize}

\item[{\rm (2)}]
For every $\varepsilon > 0$,
the following weak convergence result of deterministic measures in the weak topology of $\mathcal{M}$ holds:
$$
S^{L,M}(\nu) \stackrel{ M \to \infty }{\to } W( \sqrt{2B_L (\nu)})
$$
uniformly\footnote{We refer the reader to the original article, \cite{TsiVer}, in regards to the exact meaning of `uniform' in this sentence.}
with respect to
$(L,\nu ) \in (0, \infty ) \times \mathcal{M}$
with
$B_{L}(\nu ) \geqslant \varepsilon$.
Here,
for $\sigma > 0$,
$W(\sigma )$
stands for the Gaussian measure of which the mean function is zero
and
the variance function is given by
$$
\mathbf{E}_{W(\sigma)}[\xi (x) \xi (y)]
=\sigma^2 \mathbf{E}_W[\xi (x) \xi (y)],
$$
where $W$ is the Wiener measure.

\item[{\rm (3)}]
For every
$\varepsilon >0$,
there exist an open neighborhood
$\mathcal{U} \subset \mathcal{M}$
of $\gamma$
and
$L_{0} > 0$
with the following property:
for any
$\nu \in \mathcal{U}$
and
$L \geqslant L_0$,
it holds
that
$$
\vert
B_{L} (\nu) - B_{L} ( \gamma )
\vert
\leqslant \varepsilon .
$$
Here, $ B_{L} ( \gamma )=2 \frac{\pi -2}{\pi}$.
\end{itemize}
\end{Thm} 

\begin{Rm} 
\label{Rm:Coordinates} 
\begin{itemize}
\item[(a)]
In Section 5 of Tsirelson and Vershik (1998), the authors prove the $\mu$-nonlinearizability of $\Omega$.
Within this proof one finds the proof of the blackness property.

\item[(b)]
For future arguments, it is useful to define the following image measure
$
\mu_{k} := R_{k,\infty} (\mu)
$
which is the law of $\omega_{k}$ under $\mu$.
We will also use the law of $\xi_{k}$ defined by
the equation (\ref{Def:Xi-map}),
which will
take place in
$\nu$ in the above statements.

\item[(c)]
If we interpret that (\ref{Def:R-map})
represents a time series, one may say that
the black noise is constructed on the
remote past.
\end{itemize}
\end{Rm} 

\section{Implementation and Result}

In this section, we shall consider the black noise set-up
$
(
	\Omega (\mathbb{R}),
	\{ \mathbb{F}(s,t) \}_{s \leqslant t},
	\mu ,
	\{ T_{h} \}_{h\in\mathbb{R}}
)
$
described in Theorem~\ref{Thm:Tsi-Ver}--(1).
We will explain how to generate a
sample of sequences of continuous paths from the measure $\mu$.

\subsection{Implementation}
\label{Sec:procedure}

In order to introduce the simulation method, we define an auxiliary transformation as follows.

\begin{Def} 
\label{Def:Xi} 
We define
$\Xi_k : C(\mathbb{R}) \to C(\mathbb{R})$
by
\begin{equation*}
\begin{split}
(\Xi_k (\omega_{k+1}))(x)
:= \xi_k(x)
&= 
\displaystyle \frac{\sqrt{M-1}}{2r_{k+1}}
\int_{s_k x-r_{k+1}}^{s_k x + r_{k+1}} \omega_{k+1}(y)\mathrm{d}y,
\quad
x \in \mathbb{R}.
\end{split}
\end{equation*}
\end{Def} 

A graphical relation among
$\omega_{k}$,
$\omega_{k+1}$,
$\xi_{k}$
and
$\xi_{k+1}$
can be summarized in the diagram~(\ref{eq:diagram}).
\begin{equation}
\label{eq:diagram} 
\begin{split}
\hspace{-8mm}
\xymatrix@M=8pt{
	\text{{\large $\omega_{k}$}}
	&
	&
	&
	&
	&
	\text{{\large $\omega_{k+1}$}}
	\ar@/_/[lllll]_-{\text{$R_{k,k+1}$}}
	\ar@/_/[llllldd]_-{\text{$\Xi_{k}$}}
	\ar@/_/[llllldd]^-{\text{$
	\begin{array}{l}
	\hspace{10mm} \xi_{k} (x) \vspace{2mm}\\
	\displaystyle
	=
	\frac{ \sqrt{M-1} }{ 2 r_{k+1} }
	\int_{ s_{k} x - r_{k+1} }^{ s_{k} x + r_{k+1} }
	\omega_{k+1} (y)
	\mathrm{d} y
	\end{array}
	$}}
	\\
	\\
	\text{{\large $\xi_{k}$}}
	\ar@/^/[uu]^-{\text{$
	\begin{array}{c}
	\omega_{k} (x) \\
	\rotatebox{90}{$=$} \\
	\varphi ( L \xi_{k} ( \frac{x}{s_{k}} ) )
	\end{array}
	$}}
	&
	&
	&
	&
	&
	\text{{\large $\xi_{k+1}$}}
	\ar@/^/[lllll]^-{\text{$\displaystyle
		( R \xi_{k+1} )(x)
		:=
		\frac{1}{2\sqrt{M-1}}
		\int_{Mx-M+1}^{Mx+M-1}
		\varphi ( L\xi_{k+1}(y) )
		\mathrm{d}y
	$}}
	\ar@/_/[uu]_-{\text{$
	\begin{array}{c}
	\omega_{k+1} (x) \\
	\rotatebox{90}{$=$} \\
	\varphi ( L \xi_{k+1} ( \frac{x}{s_{k+1}} ) )
	\end{array}
	$}}
}
\end{split}
\end{equation}
Therefore, for the implementation of the black noise
$\omega = ( \omega_{k} )_{k=1}^{\infty}$,
it is enough to construct the associated $\xi_{k}$'s
because of the presence of upward arrows.
In particular, notice that the application
$
R: C(\mathbb{R}) \ni \xi_{k+1} \mapsto \xi_{k} \in C(\mathbb{R})
$
(the bottom line in the diagram~(\ref{eq:diagram}))
does not depend on $k$ when compared to that of $R_{k,k+1}$
(the top line in the diagram~(\ref{eq:diagram})).
This enables us to write codes for an implementation of $\xi_{k}$'s somewhat simply.

According to Theorem~\ref{Thm:Tsi-Ver},
the probability measure $\mu$ in Theorem~\ref{Thm:Tsi-Ver}--(1)
satisfies the following properties:
for any small open neighborhood
$
\mathcal{U} \subset \mathcal{M}
$
of $\gamma$,
if we take sufficiently large $L,M > 0$ then
\begin{itemize}
\item[$\bullet$]
Theorem~\ref{Thm:Tsi-Ver}--(1):
for $k=0,1,2, \ldots $,
\begin{equation*}
(\text{the law of $\xi_{k}$ under $\mu$})
=
\Xi_k (R_{k+1, \infty} (\mu)) \in \mathcal{U} .
\end{equation*}
\end{itemize}
Since $\mathcal{U}$ is small and $L$ is large,
we have for any $k=0, 1, 2, \ldots$,
\begin{itemize}
\item[$\bullet$]
Theorem~\ref{Thm:Tsi-Ver}--(3):
$
B_{L} ( \Xi_{k} (R_{k+1, \infty}(\mu )) )
\approx
2 \frac{\pi -2}{\pi}
$
($= B_{L} ( \gamma )$),
\end{itemize}
and hence by Theorem~\ref{Thm:Tsi-Ver}--(2),
when $M$ is large enough,
we reach that as $L,M \to \infty$, 
$
S^{L, M} ( \Xi_{k} ( R_{k+1, \infty} ( \mu ) ) )
$
is close to
$
W
(
	2 \sqrt{ \frac{\pi -2}{\pi} }
)
$
for any $ k=0, 1, 2, \ldots$.
More precisely, this can be stated as follows.

\begin{Prop} 
\label{Prop:the_contr} 

Consider all Borel probability measures on $C(\mathbb{R})$.
Let $d$ be an arbitrary distance function on them, which respects the weak convergence.
Then for any $\varepsilon > 0$, there exist positive numbers $L_{0}$ and $M_{0}$, and a probability measure $\mu$ on $\Omega (\mathbb{R})$ such that Theorem~\ref{Thm:Tsi-Ver}--(1)--(i) holds and, for any $L \geqslant L_{0}$ and $M \geqslant M_{0}$, we have
$
d \left( S^{L, M} ( \Xi_{k} ( R_{k+1, \infty} ( \mu ) ) ), W( 2 \sqrt{ \frac{\pi -2}{\pi} }) \right) < \varepsilon
$.

\end{Prop} 

\begin{proof} 
Let $\varepsilon > 0$ be an arbitrary positive real number.
Note that $(L,\nu) \mapsto W(\sqrt{2 B_{L}(\nu)})$ is continuous, and $B_{L}(\gamma) = 2 \frac{\pi -2}{\pi}$ for every $L>0$.
By Theorem~\ref{Thm:Tsi-Ver}--(3), there exist $L_{1}>0$ and an open neighborhood $\mathcal{U}_{1}$ of $\gamma$ such that
$\displaystyle
\sup_{\substack{L\geqslant L_{1},\\ \nu \in \mathcal{U}_{1}}}
d \left( W(\sqrt{2B_{L}(\nu)}), W\big( 2 \sqrt{{\textstyle \frac{\pi -2}{\pi}}} \big) \right)
<
{\textstyle \frac{\varepsilon}{2}}
$.
By Theorem~\ref{Thm:Tsi-Ver}--(2), we can take $M_{1}>0$ such that
$\displaystyle
\sup_{\substack{L,\, \nu :\\ B_{L}(\nu) \geqslant \frac{\pi -2}{\pi}}}
d \left( S^{L,M}(\nu ), W\big( \sqrt{2B_{L}(\nu)} \big) \right)
<
{\textstyle \frac{\varepsilon}{2}}
$
for any $M \geqslant M_{1}$.
By Theorem~\ref{Thm:Tsi-Ver}--(3), there exist $L_{2}>0$ and an open neighborhood $\mathcal{U}_{2}$ of $\gamma$ such that
$
\vert B_{L}(\nu) - 2 \frac{\pi -2}{\pi} \vert \leqslant \frac{\pi -2}{\pi}
$
(which implies that $B_{L}(\nu) \geqslant \frac{\pi -2}{\pi}$)
for every $L \geqslant L_{2}$ and $\nu \in \mathcal{U}_{2}$.
Finally, by Theorem~\ref{Thm:Tsi-Ver}--(1), there exist $L_{3}>0$, $M_{2}>0$ and a probability measure $\mu$ on $\Omega (\mathbb{R})$ such that
the tuple $( \Omega (\mathbb{R}) ,\{ \mathbb{F}(s,t) \}_{s \leqslant t}, \mu , \{ T_{h} \}_{h \in \mathbb{R}} )$ is a black noise and
$
\Xi_{k} ( R_{k+1,\infty} (\mu) ) \in \mathcal{U}_{1} \cap \mathcal{U}_{2}
$
for every $k$.

Now, by setting $L_{0} := \max \{ L_{1}, L_{2}, L_{3} \}$ and $M_{0} := \max \{ M_{1}, M_{2} \}$,
we have, for every $L \geqslant L_{0}$ and $M \geqslant M_{0}$,
\begin{equation*}
\begin{split}
&
d \left( S^{L,M} ( \Xi_{k} ( R_{k+1,\infty} (\mu) ) ), W \big( 2 \sqrt{{\textstyle \frac{\pi -2}{\pi} }} \big) \right) \\
&\leqslant
\sup_{\substack{L,\nu :\\ B_{L}(\nu) \geqslant \frac{\pi -2}{\pi}}}
d \left( S^{L,M}(\nu ), W ( \sqrt{2B_{L}(\nu )} ) \right)
+
\sup_{\substack{L \geqslant L_{0}, \\ \nu \in \mathcal{U}_{1}\cap \mathcal{U}_{2}}}
d \left( W(\sqrt{ 2 B_{L}(\nu) }), W \big( 2 \sqrt{{\textstyle \frac{\pi -2}{\pi} }} \big) \right) \\
&<
\frac{\varepsilon}{2} + \frac{\varepsilon}{2} = \varepsilon .
\end{split}
\end{equation*}
\end{proof} 

With this in mind,
our implementation goes as
follows:
Fix $N \in \mathbb{N}$ and consider the following approximation to diagram (\ref{eq:proj_sys})
\begin{equation*}
\Omega_1 (\mathbb{R})
\stackrel{R_{1,2}}{\leftarrow}
\Omega_2 (\mathbb{R})
\stackrel{R_{2,3}}{\leftarrow}
\cdots
\stackrel{R_{N-2,N-1}}{\leftarrow}
\Omega_{N-1} (\mathbb{R})
\stackrel{R_{N-1, N}}{\leftarrow}
\Omega_N (\mathbb{R}).
\end{equation*}
If we incorporate the diagram in (\ref{eq:diagram}), we obtain the following approximation scheme.
\begin{equation}
\label{eq:diagram2} 
\begin{split}
\hspace{-8mm}
\xymatrix@M=8pt{
	&
	\text{{\large $\omega_{1}$}}
	&
	\text{{\large $\omega_{2}$}}
	\ar@{|->}[l]_-{\text{$R_{1,2}$}}
	&
	\text{{\large $\cdots$}}
	\ar@{|->}[l]_-{\text{$R_{2,3}$}}
	&
	\text{{\large $\omega_{N-1}$}}
	\ar@{|->}[l]_-{\text{$R_{N-2,N-1}$}}
	& \text{{\large $\omega_{N}$}}
	\ar@{|->}[l]_-{\text{$R_{N-1,N}$}}
	\\
	\text{{\large $S^{L,M}\xi_{1}$}}
	&
	\text{{\large $\xi_{1}$}}
	\ar@{|->}[l]_-{\text{$S^{L,M}$}}
	\ar@{|->}[u]_-{}
	&
	\text{{\large $\xi_{2}$}}
	\ar@{|->}[l]_-{\text{$R$}}
	\ar@{|->}[u]_-{}
	&
	\text{{\large $\cdots$}}
	\ar@{|->}[l]_-{\text{$R$}}
	&
	\text{{\large $\xi_{N-1}$}}
	\ar@{|->}[l]_-{\text{$R$}}
	\ar@{|->}[u]_-{}
	&
	\text{{\large $\xi_{N}$}}
	\ar@{|->}[l]_-{\text{$R$}}
	\ar@{|->}[u]_-{}
}
\end{split}
\end{equation}
Then we build a simulation method using a Bayesian statistical model. 
For this, it will be useful to recall Remarks~\ref{Rm:CNN} and \ref{Rm:Coordinates}--(b) and keep in mind the above diagram.
The steps in the procedure are as follows:

{\bf Simulation procedure}
\begin{itemize}
	\item[(1)]
A prior distribution of $\xi_{N}$ is proposed. 
We denote it by
	$\nu_{\text{prior}}^{(N)}$.

	\item[(2)]
	From $\xi_{N}$,
	one computes in a backward way $\xi_{N-1}, \ldots , \xi_{2},\xi_{1}$
	as functions of $\xi_{N}$ using the application $R$
	in the bottom arrow of diagram (\ref{eq:diagram}).
	After this, one computes $S^{L,M}(\xi_1)$ using Definition~\ref{Def:SB}.
	Note that $S^{L,M}(\xi_1)$ is a function of $\xi_{N}$.

	\item[(3)]
	One builds the likelihood function
	$
	\mathrm{lhd} ( \xi_{N}; w )
	=
	\delta_{w} ( S^{L,M} \xi_{1} )
	$
	for $w \in C(\mathbb{R})$.
	When we actually implement this, we choose a partition in order  
	to discretize the domain of the integral in the definition of   $S^{L,M}\xi_{1}$. Then the likelihood function will be substituted by a Gaussian kernel on
	$\mathbb{R}^{d'}$,
	where $d'$ is the size of the partition of the integral domain,
	with (the discretization of) $S^{L,M} \xi_{1}$ as the mean vector and a very small covariance matrix
	since the delta functional $\delta_{w}$ can not be implemented exactly.

	\item[(4)]
	Using this construction one computes the Bayesian posterior distribution
	$$
	\nu_{\mathrm{post}}^{(N)}
	( \mathrm{d} \xi_{N} \mid w )
	:=
	\nu_{\mathrm{prior}}^{(N)}
	( \mathrm{d}\xi_{N} \mid S^{L,M}\xi_{1} = w )
	\propto \mathrm{lhd} ( \xi_{N}; w ) \, \nu_{\mathrm{prior}}^{(N)} (\mathrm{d}\xi_{N}) .
	$$
	\item[(5)]
	Finally, 
	one samples $w$'s from the distribution $W ( 2 \sqrt{ \frac{\pi -2}{\pi} } )$,
	and then sample one $\xi_{N}$ from $\nu_{\mathrm{post}}^{(N)}
	( \mathrm{d} \xi_{N} \mid w )$ for each $w$.
	Then, the totality of $\xi_{N}$'s can be regarded as samples from the distribution given by
	$$
	\hat{\nu}^{(N)} ( \mathrm{d} \xi_{N} )
	:=
	\int_{C(\mathbb{R})}
	\nu_{\text{post}}^{(N)}
	( \mathrm{d} \xi_{N} \mid w)
	W ( 2 \sqrt{{\textstyle \frac{\pi -2}{\pi}}} )
	( \mathrm{d}w ).
	$$
\end{itemize}

Then we interpret the image measure of $\hat{\nu}^{(N)}$ by
$\xi_{N} \mapsto \omega_{N}$
(the rightmost up arrow in the diagram~(\ref{eq:diagram2}))
as an approximation of $\mu$,
and the sequence $(\omega_{N}, \omega_{N-1}, \ldots , \omega_{1} )$ will be thought of as an approximating sample of the black noise.

Now, how can we recognize the blackness?
The following result due to Tsirelson (2004) would be a candidate:

\begin{Thm}[{\cite[Lemma{ ~6d9  }, Corollary{ ~6d8  }, Corollary{ ~6d11  }]{Tsi2}}]
\label{Thm:criterion} 
For any noise
$
\mathbf{N}
=
(
	\{ \mathcal{F}_{s,t} \}_{s \leqslant t},
	\{ T_{h} \}_{h \in \mathbb{R}}
)
$,
the following are equivalent:
\begin{itemize}
\item[{\rm (1)}]
$\mathbf{N}$ is black.

\item[{\rm (2)}]
For any $f \in L_{2}( \mathcal{F}_{-\infty , \infty} )$,
$$
\mathbf{H}_{1} (f)
:= 
\lim_{\{ t_1, t_2, \cdots, t_n\}\uparrow}
\sum_{l=1}^{n+1}
\mathrm{Var}
(
	\mathbf{E}
	[ f \mid \mathcal{F}_{t_{l-1}, t_{l}} ]
)
=0 ,
$$
with setting the convention: $t_{0} = -\infty$ and $t_{n+1} = +\infty$.
Here, the limit is taken over the net of all finite sets
$\{ t_{1}, t_{2}, \ldots , t_{n} \} \subset \mathbb{R}$,
$n \in \mathbb{N}$,
with $t_{1} < t_{2} < \cdots < t_{n}$ ordered by inclusion and
in such a way that the partitions tend to be dense in $\mathbb{R}$.
\end{itemize}
\end{Thm} 

Historically, apart from the theory of noises, the concepts of stability and sensitivity for a sequence of Boolean functions were originally introduced by Benjamini–Kalai–Schramm (\cite{BKS}) and derived a sufficient condition for sensitivity. 
The above result is an aspect of a beautiful intersection (\cite[Corollary~6d14]{Tsi2}) of Benjamini–Kalai–Schramm's theory of noise-sensitivity and Tsirelson's theory of noises.

Of course it is impossible to check the behaviour of $\mathbf{H}_{1}(f)$
for {\it all} $f \in L_{2}(\mathbb{F}(\mathbb{R}))$,
which is necessary for the complete verification of the blackness.
But to get some partial evidence for the blackness,
we shall consider
$
f = \int_{0}^{1} \text{\lq $\omega (x)$\rq} \mathrm{d}x
$
as an example.
This is a very fundamental example of random variables since
it stands for the increment of a Brownian motion during the time interval $[0,1]$
when we deal with a Gaussian white noise
(see Example~\ref{Eg:white_noise}).
But in this case,
since our scheme is based on the projective system
(\ref{eq:proj_sys}),
$
\omega = (\omega_1, \omega_2, \omega_3 , \ldots ) \in \Omega (\mathbb{R})
$
is not a function and hence
$\omega (x)$
does not make any sense as a real number.
Instead we shall consider
$
f_{k} = \int_{0}^{1} \omega_{k}(x) \mathrm{d}x
$,
$k=1,2,
\ldots
$, and
see how $\mathbf{H}_{1} ( f_{k} )$ behaves as $k$ grows.

\begin{Eg}[Simulation results for a white noise] 

If $\omega$ is a white noise
then by It\^{o}'s isometry we have
$
\mathbf{H}_{1} ( \int_{0}^{1} \omega (x) \mathrm{d}x )
=
\int_{0}^{1} \mathrm{d}x = 1
$.
In \cite[Section~5]{TsiVer},
Tsirelson and Vershik
describe
a white noise
as the projective limit of a projective system
(\ref{eq:proj_sys}),
where
$
R_{k-1,k}:
\Omega_{k}(\mathbb{R}) \ni \omega_{k}
\mapsto
\omega_{k-1} \in \Omega_{k-1}(\mathbb{R})
$
is replaced by
\begin{equation*}
\omega_{k-1}(x)
=
\int_{x-r_{k}}^{x+r_{k}}
V_{k} (y-x)
\omega_{k} (y)
\mathrm{d}y
\end{equation*}
with using a nonnegative function
$V_{k}$
concentrated on $(-r_{k}, r_{k})$
satisfying
$
\int V_{k} (x) \mathrm{d}x = 1
$.

In order to simulate these quantities, let
$M>1$
and
we replace $r_{k}$ by $r_{k} = M^{-k}$,
which is one slightly modified from that appeared in (\ref{Def:R-map}).
We take the nonnegative function
$V_{k}:\mathbb{R} \to [0,\infty )$
as
\begin{equation*}
V_{k} (x)
=
\left\{\begin{array}{ll}
- r_{k} \vert x \vert + r_{k}^{2}  & \text{if $\vert x \vert < r_{k}$,} \\
0 & \text{if $\vert x \vert \geqslant r_{k}$.}
\end{array}\right.
\end{equation*}
(See Figure~\ref{graph(V)}.)

\begin{figure}[H]
\begin{minipage}{0.49\hsize}
\begin{center}
\includegraphics[width=8cm]{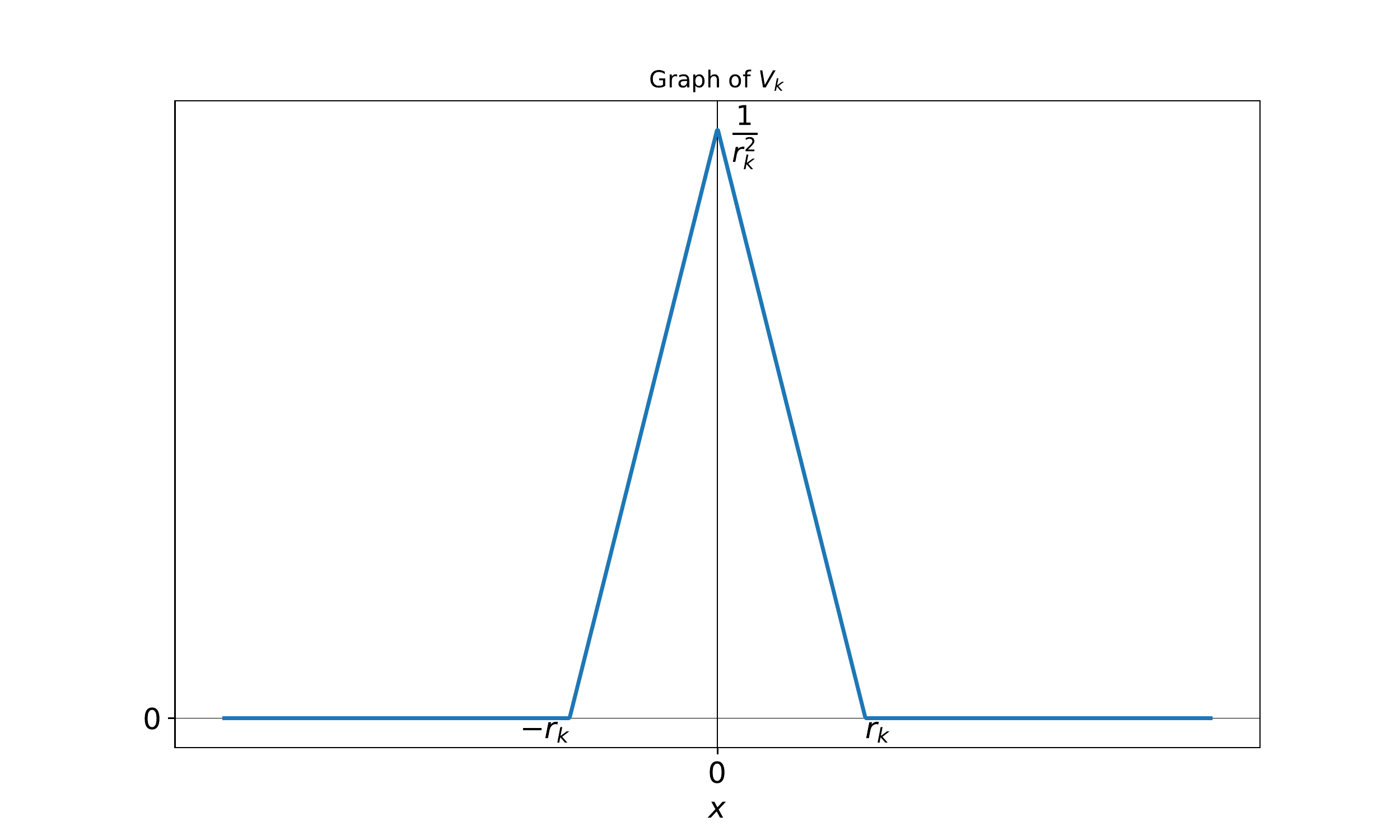}
\caption{Graph of $V_{k}$.}
\label{graph(V)}
\end{center}
\end{minipage}
\begin{minipage}{0.49\hsize}
\begin{center}
\includegraphics[width=8cm]{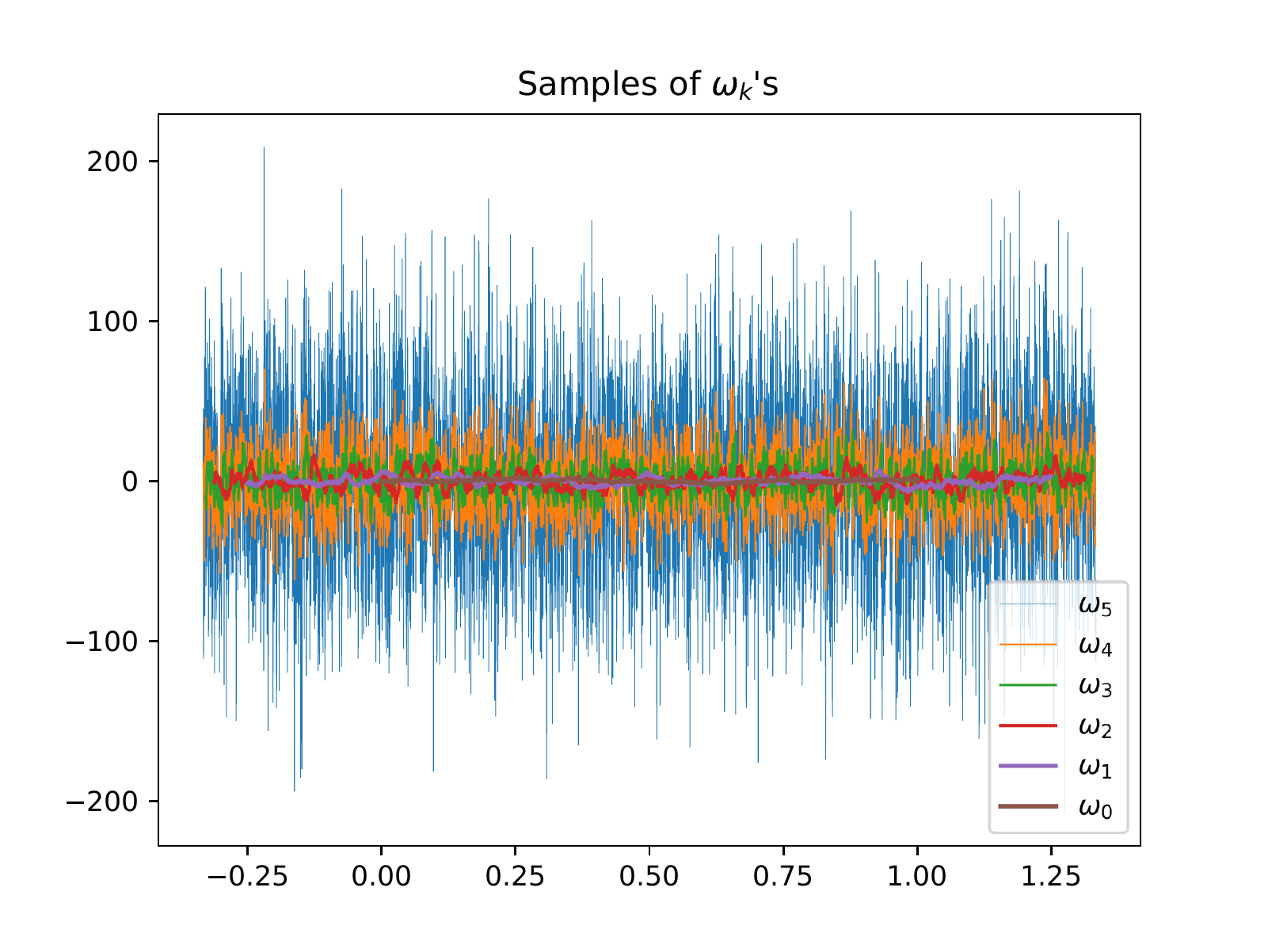}
\caption{Graphs of samples of each $\omega_{k}$'s.}
\label{graph(WN)}
\end{center}
\end{minipage}
\end{figure}
\begin{figure}[H]
\begin{center}
\includegraphics[width=10cm]{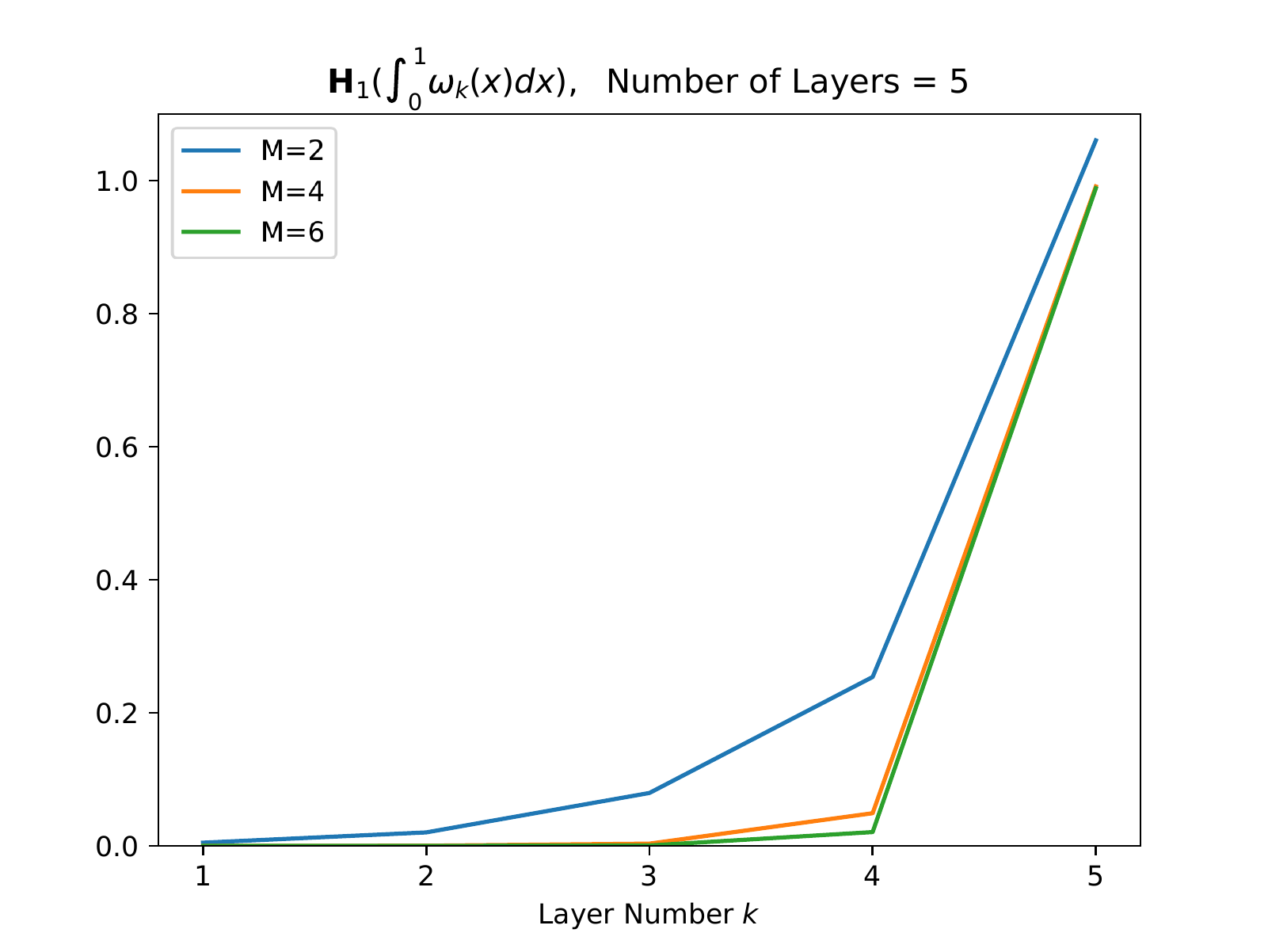}
\caption{Graphs of $k \mapsto \mathbf{H}_{1} ( \int_{0}^{1}\omega_{k}(x) \mathrm{d}x)$ for several values of $M$.}
\label{H1(WN)}
\end{center}
\end{figure}

In Figure~\ref{H1(WN)}, one sees that
$\mathbf{H}_{1}( \int_{0}^{1} \omega_{k}(x) \mathrm{d}x )$
approaches $1$ as $k$ grows.
This shows the theoretical fact that
$\mathbf{H}_{1}( \int_{0}^{1} \omega (x) \mathrm{d}x ) = 1$.
This finishes the simulation for white noise.
\end{Eg} 

In general,
$\mathbf{H}_{1}$
measures the squared $L_{2}$-norm
of a random variable in
$
L_{2} ( \mathcal{F}_{-\infty , \infty} )
$
projected onto the space of first chaos.
This means, in the situation that
$\mathbf{H}_{1} (f) > 0$ for some
$
f \in L_{2} ( \mathcal{F}_{-\infty , \infty} )
$,
one may construct a white noise as in the above example or a classical noise in Example~\ref{Eg:Levy}.

\subsection{Results}
In this section, we perform the actual simulation of one approximation to $ \omega  \in \Omega (\mathbb{R})$ which is sampled from the measure $\mu$.

We take the activation function
$\varphi:\mathbb{R} \to \mathbb{R}$
as
\begin{equation*}
\varphi (x)
=
\left\{\begin{array}{ll}
-1 & \text{if $x \leqslant -1$,} \\
x & \text{if $-1 < x < 1$,} \\
1 & \text{if $1 \leqslant x$,}
\end{array}\right.
\end{equation*}
which satisfies the
conditions~(b)--(i, ii, iii) in Section~\ref{Sec:Construction}.
The following Figure~\ref{Graph(phi)} is a graph of $\varphi$.

\begin{figure}[h]
\begin{center}
\includegraphics[width=11cm]{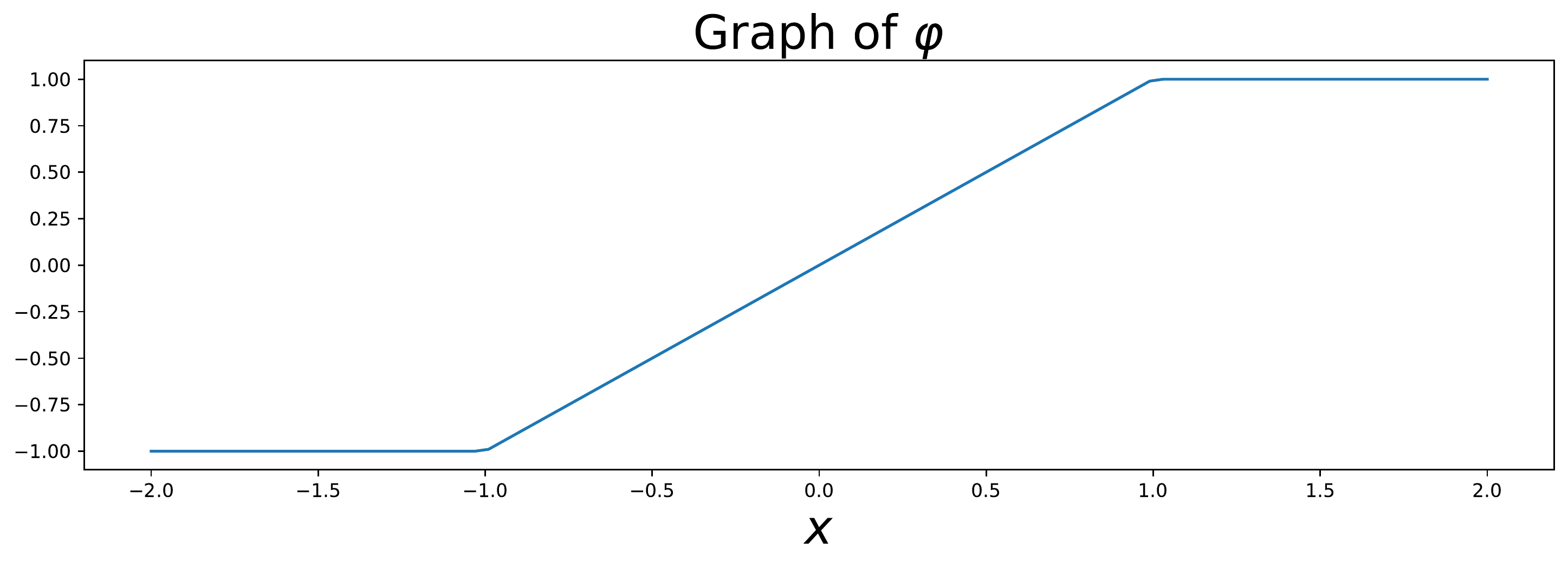}
\caption{Graph of $\varphi$}
\label{Graph(phi)}
\end{center}
\end{figure}
\begin{figure}[H]
\begin{center}
\includegraphics[width=16cm]{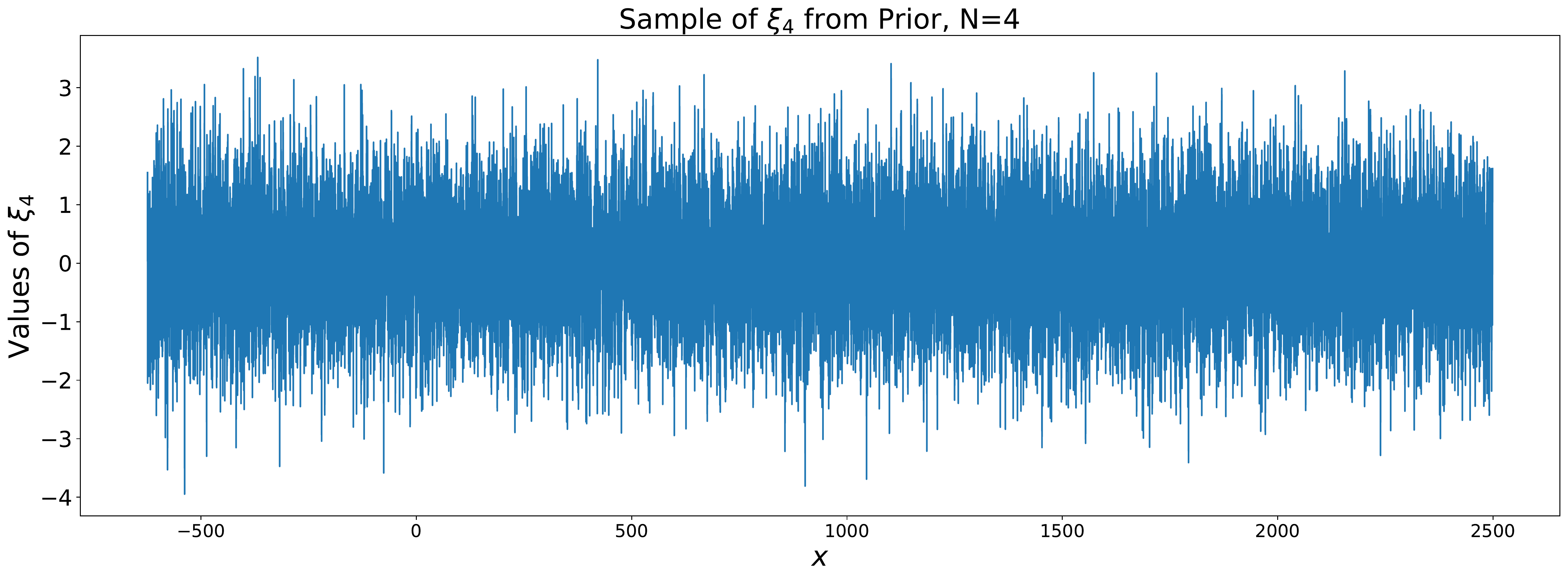}
\caption{A sample point (function) from the prior distribution on $\xi_{N}$}
\label{Prior}
\end{center}
\end{figure}

Let $M > 1$ and $L = \sqrt{M-1}$.
Firstly, we set $N = (\text{the number of layers}) = 4$.
We put a multivariate normal distribution on $\xi_{N}$
as its prior distribution.
More precisely, we have taken
$\mathrm{N}(0,1)^{\otimes d}$, here $d$ is the number of partition points in the discretization of the integral $ \int_{Mx-M+1}^{Mx+M-1} \varphi ( L\xi_{N}(y) ) \mathrm{d}y $ which appears in the definition of $\xi_{N-1}$
and $\otimes$ denotes the product of measures.
We display a sample point (function) of $\xi_{N}$ in Figure~\ref{Prior}.
For latter calculation of $\mathbf{H}_{1} (f_{k})$'s,
we need values of $\omega_{k}$'s on the interval $[0,1]$.
To obtain the values of $\omega_{1}(x)$ for $0 \leqslant x \leqslant 1$,
we needed the values of $\xi_{N}$ on $[-M^{N}+1, M^{N}(M-1)-1]$.
In Figure~\ref{Prior}, we took the case $N=4$, $M=5$ and $d = M(M^{N+1}-1)$,
and the graph of $\xi_{N}$ on the interval
$
[-M^{N}+1, M^{N}(M-1)-1]
=[ -624, 2499 ]
$ is displayed.

\begin{figure}[H]
\begin{center}
\includegraphics[width=16cm]{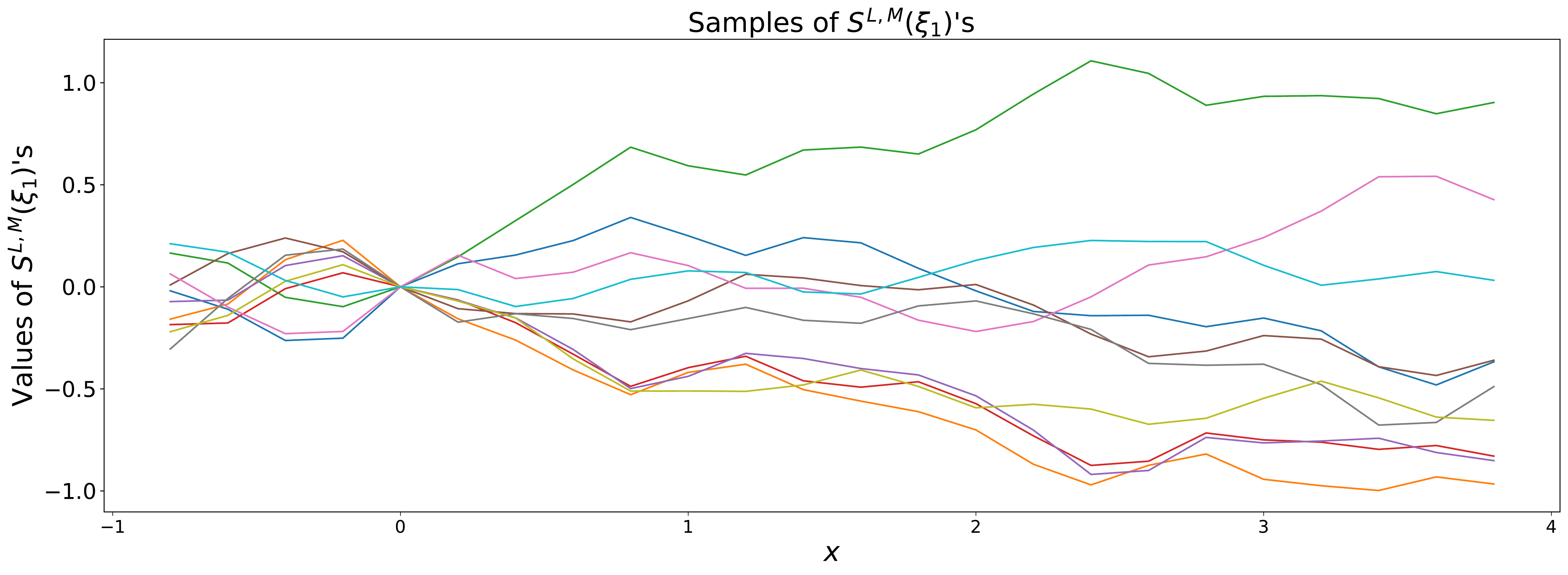}
\caption{Sample points (functions) of $S^{L,M}(\xi_{1})$
under the prior distribution of $\xi_{N}$.}
\label{S^LM}
\end{center}
\end{figure}

For a given sample point (function) of $\xi_{N}$, 
we compute the corresponding $S^{L,M} ( \xi_{1} )$.
To recognize the distribution of $S^{L,M} ( \xi_{1} )$ under the prior distribution on $\xi_{N}$, we take some samples of $\xi_{N}$. 
Figure~\ref{S^LM} is the case when $M=5$.
Note that once the values of $\xi_{N}$ on the interval $[-M^{N}+1, M^{N}(M-1)-1]$,
we can compute the values of $S^{L,M}(\xi_{1})$
on
$
[ -1+\frac{1}{M}, (M-1)-\frac{1}{M} ]
=
[-0.8, 3.8]
$
as explained in step (2) of the simulation procedure.

\begin{figure}[H]
\begin{center}
\includegraphics[width=16cm]{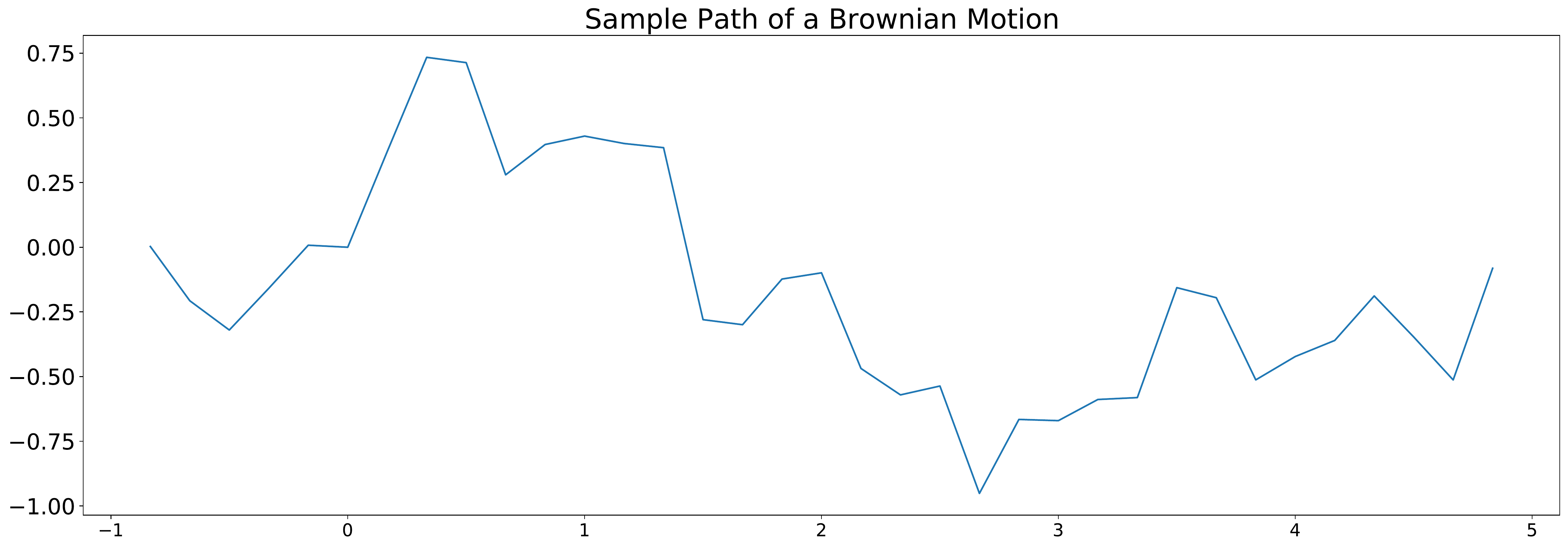}
\caption{A sample path $w$ of a Brownian motion used for Bayesian update of $\xi_{N}$ in step (4).}
\label{BM}
\end{center}
\end{figure}

\begin{figure}[H]
\begin{center}
\includegraphics[width=16cm]{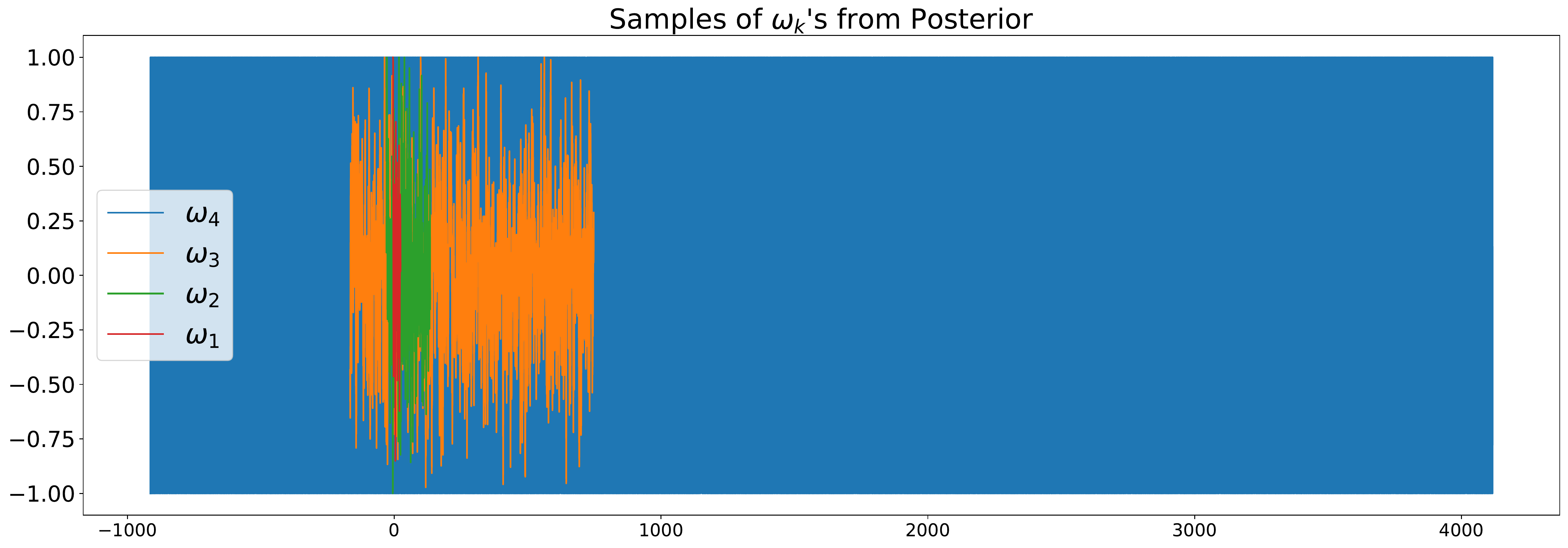}
\caption{Sample of $\omega_{N}, \omega_{N-1},\ldots , \omega_{1}$ under the posterior distribution of $\xi_{N}$}
\label{Post}
\end{center}
\end{figure}

After $S^{L,M}\xi_{1}$ is implemented as a function of $\xi_{N}$,
we can build the likelihood function as in step (3),
and then for each (discretized) $w \in C(\mathbb{R})$,
using, e.g., \texttt{PyMC} (\cite{PyMC}),
we can take a sample from the posterior distribution of $\xi_{N}$ which is defined in the step (4).
Figure~\ref{Post} is the graphs of a sample of 
$\omega_N, \omega_{N-1}, \ldots ,\omega_{1}$,
which are computed according to the diagram (\ref{eq:diagram})
under the posterior distribution of $\xi_{N}$.
Here, a sample path $w$ of a Brownian motion (scaled by multiplying $2 \sqrt{\frac{\pi - 2}{\pi}}$),
which is described in Figure~\ref{BM},
is used to construct the posterior distribution of $\xi_{N}$.


These $\omega_{k}$'s are calculated from $\xi_{k}$'s
according to the up arrows in the diagram~(\ref{eq:diagram})
and the value $\xi_{k}(x)$ at a point $x$ is calculated
by using values of $\xi_{k+1}$ on the interval $[Mx+M-1, Mx-M+1]$.
This explains the reason why the definition domain of $\omega_{k}$ becomes wider
than that of $\omega_{k-1}$.

\begin{figure}[H]
\begin{center}
\includegraphics[width=10cm]{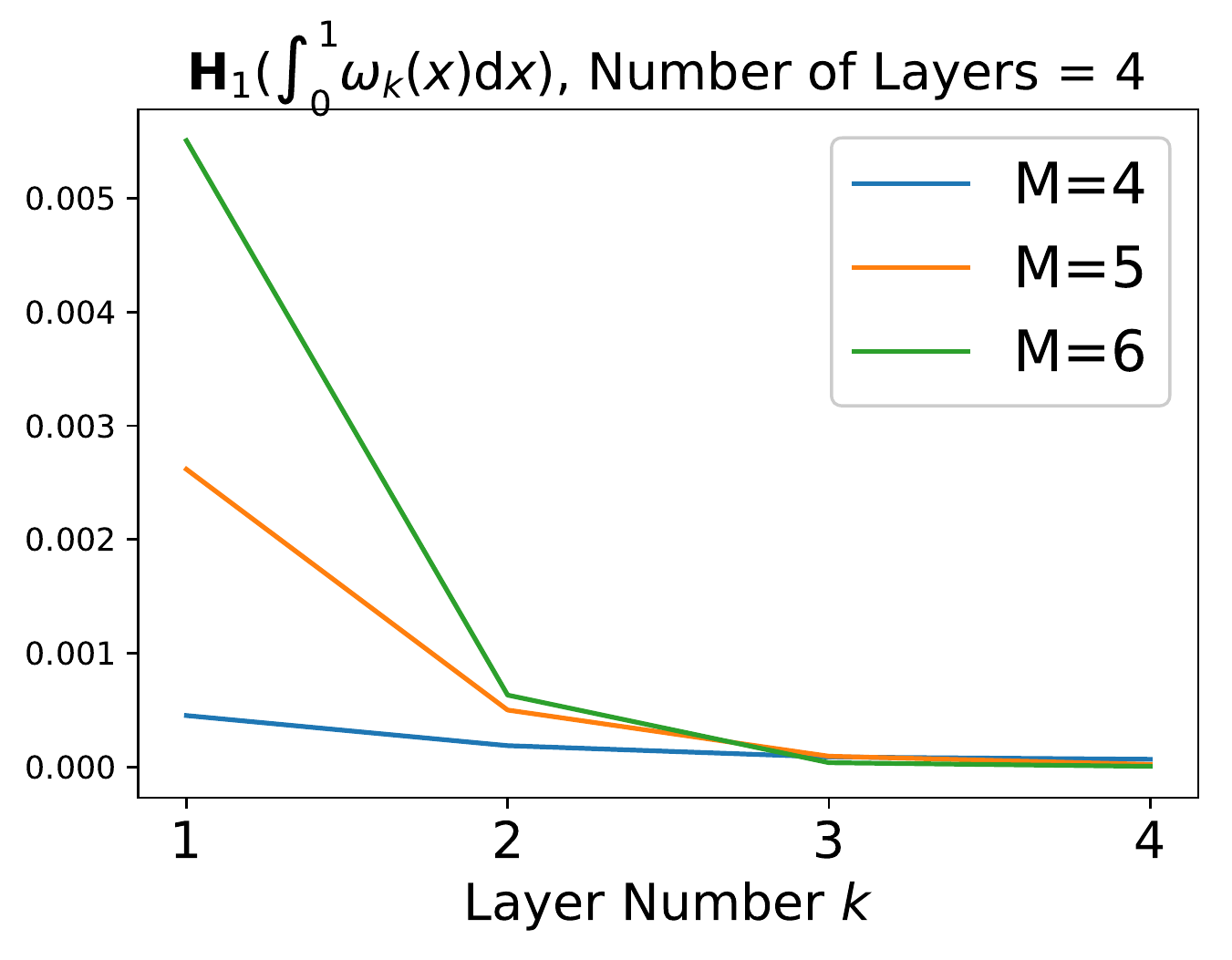}
\caption{Graphs of $k \mapsto \mathbf{H}_{1} ( \int_{0}^{1} \omega_{k} (x) \mathrm{d}x  )$ in the case of $N=4$ for several values of $M$}
\label{N=4}
\end{center}
\end{figure}

The procedure of Bayesian update above gives one sample of $\omega_{k}$ for each $k$.
We then do the procedure explained in step (5) of the simulation method.
Namely, we repeatedly sample $w$ from $W(2\sqrt{\frac{\pi -2}{\pi}})$
and pass the above procedures (1)--(4) with each $w$ to obtain many samples of $\omega_{k}$'s,
which is used for point estimations of $\mathbf{H}_{1}(f_{k})$'s described in Figure~\ref{N=4}.
Here, for the computation of $\mathbf{H}_{1}(f_{k})$'s,
we need to compute the conditional expectation
$\mathbf{E}_{\mu}[ \int_{0}^{1} \omega_{k}(x) \mathrm{d}x \mid \mathbb{F}(t_{l-1},t_{l}) ]$
(the symbol $\mathbf{E}_{\mu}$ denotes the expectation with respect to the probability measure $\mu$).
For this, we shall note the following decomposition
($s_{k}$ is the number defined in Section~\ref{Sec:Construction})
\begin{equation*}
\begin{split}
\int_{0}^{1} \omega_{k} (x) \mathrm{d}x
&=
\underbrace{
	\int_{0}^{t_{l-1}-2s_{k}} \omega_{k} (x) \mathrm{d}x
}_{
	\text{independent of $\mathbb{F} ( t_{l-1}, t_{l} )$}
}
+
\underbrace{
	\int_{ t_{l-1}-2s_{k} }^{ t_{l-1}+s_{k} } \omega_{k} (x) \mathrm{d}x
}_{
	\text{close to $0$ when $k$ is large}
} \\
&+
\underbrace{
	\int_{ t_{l-1}+s_{k} }^{ t_{l}-s_{k} } \omega_{k} (x) \mathrm{d}x
}_{
	\text{$\mathbb{F} ( t_{l-1}, t_{l} )$-measurable}
}
+
\underbrace{
	\int_{ t_{l}-s_{k} }^{ t_{l}+2s_{k} } \omega_{k} (x) \mathrm{d}x
}_{
	\text{close to $0$ when $k$ is large}
}
+
\underbrace{
	\int_{ t_{l}+2s_{k} }^{ 1 } \omega_{k} (x) \mathrm{d}x
}_{
	\text{independent of $\mathbb{F} ( t_{l-1}, t_{l} )$}
},
\end{split}
\end{equation*}
from which it is natural to estimate
$\mathbf{E}_{\mu}[ \int_{0}^{1} \omega_{k}(x) \mathrm{d}x \mid \mathbb{F}(t_{l-1},t_{l}) ]$
by
$\int_{ t_{l-1}+s_{k} }^{ t_{l}-s_{k} } \omega_{k} (x) \mathrm{d}x$.
Then the point estimation of
$\mathbf{H}_{1}(\int_{0}^{1} \omega_{k} (x) \mathrm{d}x)$
is done using samples from the posterior distribution of $\omega_{k}$
obtained by Monte Carlo methods.

\noindent
One may see that
the value of $\mathbf{H}_{1} ( \int_{0}^{1} \omega_{k} (x) \mathrm{d}x  )$ approaches $0$
as $k$ grows, which would be a numerical evidence
of the validity of Theorem~\ref{Thm:criterion}--(2)
in the case of $f = \int_{0}^{1} \text{`}\omega (x)\text{'} \mathrm{d}x$.

\begin{figure}[H]
\begin{center}
\includegraphics[width=10cm]{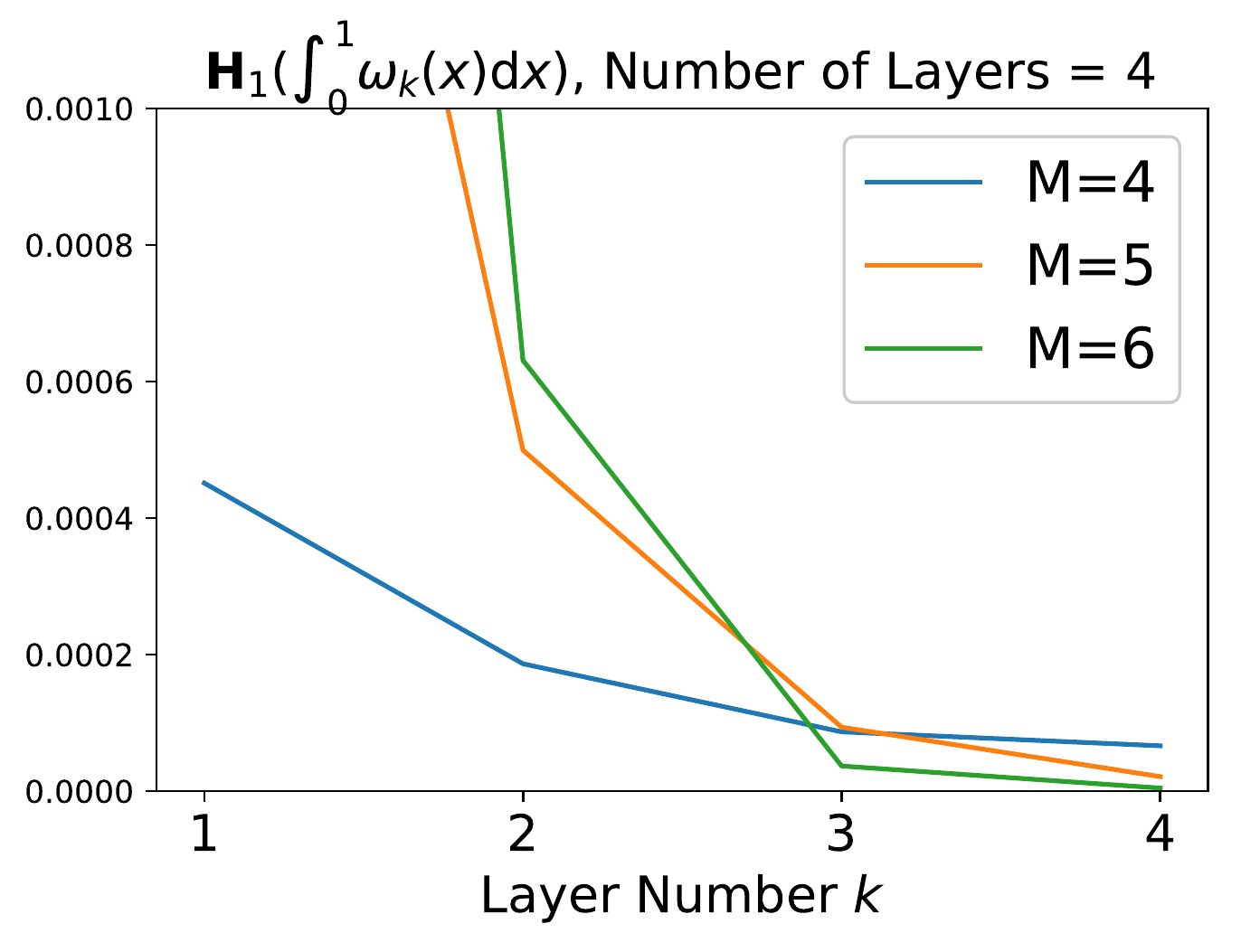}
\caption{Zoomed version of Figure~\ref{N=4}}
\label{N=4v2}
\end{center}
\end{figure}

Figure~\ref{N=4v2} focuses on the
behavior of $ \mathbf{H}_{1} ( \int_{0}^{1} \omega_{k} (x) \mathrm{d}x  ) $
near zero.
\noindent
In particular, the larger the value of $M$,
the smaller the value of $\mathbf{H}_{1} ( f_{4} )$
which would imply the darker the noise.

From the above experiments, we conclude that we have obtained circumstantial evidence to claim we have simulated the black noise constructed by Tsirelson and Vershik (1998).
Also, Figure~\ref{Post} visually confirms that Tsirelson–Vershik's black noise has a similar property of oscillating infinitely fast in space, like the paradoxical motion of an ideal incompressible fluid described by Shnirelman quoted by Tsirelson, which is also quoted in the introduction.

Finally, we show run-time to obtain a single sample $\omega_{N}$ by performing a Bayesian update. 
We show them separately by specific values of parameters. 
Note that given the value of $M$, we had determined that $L=\sqrt{M-1}$.
All calculations were performed by using Google Colaboratory (\texttt{Colab Pro+}).

\begin{table}[h]
 \caption{Run-time to obtain one sample $\omega_{N}$}
 \label{table:runtime}
 \centering
  \begin{tabular}{ccc}
   \hline
   $N$ & $M$ & Run-Time \\
   \hline \hline
   3 & 4.0 & 00h05m55s \\
   3 & 5.0 & 00h46m59s \\
   3 & 6.0 & 03h11m22s \\
   \hline
   4 & 4.0 & 01h41m50s \\
   4 & 5.0 & 18h32m15s\\
   4 & 6.0 & Crashed (Shortage of RAM) \\
   \hline
  \end{tabular}
\end{table}
The numerical results in Figures~\ref{N=4} and \ref{N=4v2} for $M=6$ shown above were run on a different machine with more memory. Here, we showed the results on \texttt{Colab Pro+} to compare execution speeds in a more standard environment that anyone can use.

\section*{Acknowledgments}
The authors would like to express their sincere appreciations to Professor Arturo Kohatsu-Higa for his valuable comments.


\begin{thebibliography}{00}

\bibitem{AIW}
{\sc J.~Akahori}, {\sc M.~Izumi}, {\sc S.~Watanabe},
``Noises, stochastic flows and $E_{0}$-semigroups,"
in:
Selected papers on probability and statistics
(Papers translated from the Japanese, originally published in S\={u}gaku),
Amer. Math. Soc. Transl. Ser. 2, {\bf 227},
Amer. Math. Soc., Providence, RI,
2009, pp.~1--23.
{\tt http://dx.doi.org/10.1090/trans2/227}.


\bibitem{Ar}
{\sc N.~Aronszajn},
``Theory of Reproducing Kernels,"
Trans. Amer. Math. Soc.
{\bf 68}(3)
(1950),
pp.~337--404.
{\tt doi:10.1090/S0002-9947-1950-0051437-7}%


\bibitem{BKS}
{\sc I.~Benjamini}, {\sc G.~Kalai} and {\sc O.~Schramm},
``Noise sensitivity of Boolean functions and applications to percolation,"
Inst. Hautes \'{E}tudes Sci. Publ. Math. no. 90 (1999), pp.~5--43.


\bibitem{ElFe}
{\sc T.~Ellis} and {\sc O.~N.~Feldheim},
``The Brownian web is a two-dimensional black noise,"
Annales de l'Institut Henri Poincar\'{e}, Probabilit\'{e}s et Statistiques.
{\bf 52}(1) (2016), pp.~162--172.


\bibitem{Ito1}
{\sc K.~It\^{o}},
``Stochastic Integral," Proc. Imperial Acad. Tokyo {\bf 20} (1944), pp.~519--524.%


\bibitem{Ito2}
{\sc K.~It\^{o}},
``On stochastic differential equations," Memoirs, American Mathematical Society {\bf 4} (1951), pp.~1--51.%


\bibitem{Ka1}
{\sc G.~Kalai},
``Is the Universe Noise-Sensitive?,"
arXiv preprint
(2007).
\texttt{arXiv:hep-th/0703092}


\bibitem{Ka2}
{\sc G.~Kalai},
``Noise Sensitivity--The case of Percolation,"
presented in Hebrew University HEP seminar, 25 April 2007,
the power point presentation available from
\texttt{https://gilkalai.wordpress.com/2009/03/06/noise-sensitivity-lecture-and-tales/},
(accessed August 9th, 2022).


\bibitem{KuWa}
{\sc H.~Kunita} and {\sc S.~Watanabe},
``On square-integrable martingales,"
Nagoya Math. J. {\bf 30} (1967), pp.~209--245.%


\bibitem{LeRa1}
{\sc Y.~Le~Jan} and {\sc O.~Raimond},
``Flows, coalescence and noise,"
The Annals of Probability {\bf 32}(2) (2004), pp.~1247--1315.
\texttt{DOI: 10.1214/009117904000000207}%


\bibitem{LeRa2}
{\sc Y.~Le~Jan} and {\sc O.~Raimond},
``Sticky flows on the circle and their noises,"
Probability Theory and Related Fields {\bf 129}(1) (2004), pp.63--82.
\texttt{https://doi.org/10.1007/s00440-003-0324-9}%


\bibitem{LBNSPS}
{\sc J.~Lee}, {\sc Y.~Bahri}, {\sc R.~Novak}, {\sc S.~Schoenholz}, {\sc J.~Pennington}, {\sc J.~Sohldickstein},
``Deep neural networks as gaussian processes,"
published as a conference paper at In International Conference on Learning Representations
(ICLR)
2018.
(also {\tt arXiv:1711.00165v3})


\bibitem{Ne}
{\sc R.~M.~Neal},
``Priors for infinite networks,"
in: Bayesian Learning for Neural Networks,
Lecture Notes in Statistics, Vol. 118
(Springer, New York, 1996),
1st ed.,
pp.~29--53.


\bibitem{Oga1}
{\sc S.~Ogawa},
``Quelques propri\'{e}t\'{e}s de l'int\'{e}grale stochastique du type noncausal,"
Japan J. Appl. Math., {\bf 1} (1984), pp.~405--416.%


\bibitem{Oga2}
{\sc S.~Ogawa},
``The stochastic integral of noncausal type as an extension of the symmetric integrals,"
Japan J. Appl. Math. {\bf 2} (1985), pp.~229--240.
\texttt{https://doi.org/10.1007/BF03167046}%


\bibitem{PyMC}
{\sc J.~Salvatier}, {\sc T.V.~Wiecki} and {\sc C.~Fonnesbeck},
``Probabilistic programming in Python using PyMC3,"
(2016),
{\tt DOI: 10.7717/peerj-cs.55}


\bibitem{SSG}
{\sc O.~Schramm}, {\sc S.~Smirnov} and {\sc C.~Garban},
``On the scaling limits of planar percolation,"
Ann. Probab. {\bf 39}(5) (2011), pp.~1768--1814.
\texttt{DOI: 10.1214/11-AOP659}


\bibitem{Shn}
{\sc A.~Shnirelman},
``On the nonuniqueness of weak solution of the Euler equation,"
Comm. Pure Appl. Math., {\bf 50}(12) (1997), pp.~1261--1286.


\bibitem{TsiVer}
{\sc B.~Tsirelson}, {\sc A.~M.~Vershik},
``Examples of nonlinear continuous tensor product of measure spaces and non-Fock factorizations,"
Rev. Math. Phys.
{\bf 10}(1)
(1998),
pp.~81--145.
{\tt https://doi.org/10.1142/S0129055X98000045}

\bibitem{Tsi1}
{\sc B.~Tsirelson},
``Unitary Brownian motions are linearizable,"
arXiv preprint
(1998).
\texttt{arXiv:math.PR/9806112}
(also MSRI Preprint No.~1998-027)


\bibitem{Tsi2}
{\sc B.~Tsirelson},
``Scaling limit, noise, stability,"
\texttt{arXiv:math/0301237},
in:
J.~Picard (Ed.),
Lectures on probability theory and statistics,
Lecture Notes in Math., 1840, Springer, Berlin,
2004,
pp.~1--106.
[Note: the numbers of sections, propositions, lemmas and theorems refer to the arXiv version.]


\bibitem{Tsi3}
{\sc B.~Tsirelson},
``Nonclassical stochastic flows and continuous products,"
Probability Surveys
{\bf 1}
(2004), pp.~173--298.


\bibitem{WaWa}
{\sc J.~Warren} and {\sc S.~Watanabe},
``On Spectra of Noises Associated with Harris flows,"
Stochastic analysis and related topics in Kyoto, Adv. Stud. Pure Math. {\bf 41} (2004), Math. Soc. Japan, Tokyo, pp.~351--373.


\bibitem{Wa}
{\sc S.~Watanabe},
``Analysis of Wiener Functionals (Malliavin Calculus) and its Applications to Heat Kernels,"
Annals of Probability, {\bf 30} (1987), pp~1--39.
\texttt{doi:10.1214/aop/1176992255}%


\bibitem{Wa2}
{\sc S.~Watanabe},
``A simple example of black noise,"
Bull. Sci. Math. {\bf 125}(6–7) (2001), pp.~605--622.


\end{thebibliography}
\end{document}